\documentclass{amsart}
\usepackage{amssymb,latexsym,graphicx}

\DeclareMathOperator{\id}{id}
\DeclareMathOperator{\supp}{supp}

\def\L{\mathcal{L}}
\def\J{\mathcal{J}}
\def\P{\mathcal{P}}

\def\R{\mathbb R}
\def\Z{\mathbb Z}
\def\N{\mathbb N}
\def\C{\mathbb C}
\def\id{{\rm id}}
\def\D{{\mathcal D}}
\def\orbit{\hbox{orbit}}
\def\zeroOne{zero-one }
\def\Q{{\mathbb Q}}
\def\pwm{piecewise monotonic}
\def\htau{{\widehat\tau}}
\def\wJ{\widehat J}

\def\mapright#1{\smash{\mathop{\longrightarrow}\limits^{#1}}}
\def\mapdown#1{\downarrow\rlap{$\vcenter{\hbox{$\scriptstyle#1$}}$}}

\def\mapright#1{\smash{\mathop{\longrightarrow}\limits^{#1}}}

\theoremstyle{plain}

\newtheorem{theorem}{Theorem}[section]
\newtheorem{proposition}[theorem]{Proposition}
\newtheorem{lemma}[theorem]{Lemma}
\newtheorem{corollary}[theorem]{Corollary}

\theoremstyle{definition}
\newtheorem{definition}[theorem]{Definition}
\newtheorem{example}[theorem]{Example}

\theoremstyle{remark}
\newtheorem*{notation}{Notation}

\numberwithin{equation}{section}

\newif\ifproofing
\proofingfalse
\newcommand{\note}[1]{\ifproofing\marginpar{#1}\fi}

\newcommand{\prop}[2]{\begin{proposition}\label{#1}\note{#1}#2\end{proposition}}
\newcommand{\lem}[2]{\begin{lemma}\label{#1}\note{#1}#2\end{lemma}}
\newcommand{\cor}[2]{\begin{corollary}\label{#1}\note{#1}#2\end{corollary}}
\newcommand{\theo}[2]{\begin{theorem}\label{#1}\note{#1}#2\end{theorem}}
\newcommand{\defi}[2]{\begin{definition}\label{#1}\note{#1}#2\end{definition}}
\newcommand{\exe}[2]{\begin{example}\label{#1}\note{#1}#2\end{example}}
\newcommand{\notate}[1]{\begin{notation}#1\end{notation}}

\newcommand{\prooff}[1]{\begin{proof}#1\end{proof}}

\begin{document}

\title{Dimension groups for interval maps}
\author{Fred Shultz}
\address{Wellesley College\\Wellesley, Massachusetts 02481}
\email{fshultz@wellesley.edu}
\subjclass[2000]{Primary 37E05, 46L80}
\keywords{dimension group, interval map, piecewise monotonic, unimodal map,
tent map, $\beta$-shift, interval exchange map, C*-algebra}

\begin{abstract}

With each piecewise monotonic map $\tau$  of the unit interval,
a dimension triple is associated.  The dimension triple, viewed as a
$\Z[t,t^{-1}]$ module, is finitely generated, and generators are identified.
Dimension groups are computed for Markov maps, unimodal maps, multimodal maps,
 and  interval exchange maps.
It is shown that the dimension group defined here is isomorphic to $K_0(A)$,
where $A$ is  a C*-algebra (an ``AI-algebra") defined in dynamical terms.

\end{abstract}

\maketitle

\section{Introduction}

Given a \pwm\ map $\tau$ of the unit interval into itself,
our goal is to associate a dimension group $DG(\tau)$, providing
an  invariant for the original map.  

In a dynamical context, dimension
groups were introduced by Krieger  \cite{KrDim}. Motivated by Elliott's classification of
AF-algebras by their dimension groups \cite{Ell}, he gave a purely dynamical definition of the dimension group, via an
equivalence relation defined by the action of an ``ample" group of homeomorphisms on the space of closed and open subsets of a
zero dimensional metric space. For a shift of finite type, he associated  an ample group, and in \cite{KrShift},  showed that
the dimension triple (consisting of an ordered dimension group with a canonical automorphism) completely determines an
irreducible aperiodic shift of finite type up to shift equivalence.

Krieger's definition of a dimension group was extended
to the context of a surjective local homeo\-morphism $\sigma:X\to X$, where $X$ is a compact zero
dimensional metric space, by Boyle, Fiebig, and Fiebig (\cite{BFF}), who defined a dimension group
called the  ``images group", and used this group as a tool in  studying commuting local homeo\-morphisms.

We sketch the construction of the dimension group in \cite{BFF}.  If $\sigma:X\to X$,  the
transfer map
$\L:C(X,\Z)\to C(X,\Z)$ is defined by
\begin{equation}
\L_\sigma f(x) = \sum_{\sigma y = x} f(y).
\end{equation}
Then $G_\sigma$ is defined to be the set of equivalence classes of functions in $C(X,\Z)$, where
$f
\sim g$ if
$\L^n f = \L^n g$ for some $n \ge 0$. Addition is given by $[f]+[g] = [f+g]$, and the positive
cone consists of classes $[f]$ such that $\L^n f\ge 0$ for some $n \ge 0$. This is the same as
defining clopen sets
$E, F$ in
$X$ to be equivalent if for some $n \ge 0$, $\sigma^n$ is 1-1 on $E$ and $F$, and $\sigma^n(E) =
\sigma^n (F)$, and then building a group out of these equivalence classes by defining $[E] + [F]
= [E\cup F]$ when $E, F$ are disjoint.  This latter definition is the one given in \cite{BFF};
the authors then observe that it is equivalent to the definition above involving the transfer
operator. If $\sigma$ is surjective, this dimension group is the same as the stationary inductive limit of
$\L_\sigma:C(X,\Z) \to C(X,\Z)$, cf. Lemma \ref{1.17} and equation (\ref{(1.8)}) in the current
paper.

A related approach was taken by Renault in
\cite{Ren}.  If $X_1, X_2, \ldots$ are compact metric spaces, and $T_n:C(X_n) \to
C(X_{n+1})$ is a sequence of positive maps, Renault defined the associated dimension group to be
the inductive limit of this sequence. If $X_n = X$ for all $n$, and  $T =
\L_\sigma$ for a surjective local homeo\-morphism
$\sigma$, the resulting dimension group is formally  similar to that in \cite{BFF}. However, the
use of $C(X)$ instead of $C(X,\Z)$ results in different dimension groups.

These definitions of dimension groups aren't  directly applicable to interval maps, since 
such maps are rarely local homeo\-morphisms. We therefore associate with each piecewise monotonic map
$\tau:[0,1]\to [0,1]$ a local homeo\-morphism. This is done by disconnecting the
unit interval at a countable set of points, yielding a   space
$X$, and then lifting $\tau$ to a local homeo\-morphism $\sigma:X\to X$. The properties of $\sigma$
are closely related to those of $\tau$.  A similar technique has long been used in studying
interval maps, e.g., cf. \cite{KeaneDisc, HofDecomp, RuBook,WalBeta}.

If $\tau:I\to I$ is \pwm, and $\sigma:X\to X$ is the associated local homeo\-morphism, the space
$X$ will be a compact subset of
$\R$, but will not necessarily be zero dimensional. However the map $\sigma$ will have the 
property that each point has a clopen neighborhood on which $\sigma$ is injective, and this 
allows us to define a dimension group in  the same way as in \cite{BFF}. We then define
$DG(\tau)$ to be the dimension group $G_\sigma$ associated with $\sigma$.

 The result is a triple
$(DG(\tau), DG(\tau)^+,\L_*)$, where $DG(\tau)$ is a dimension group  with positive cone
$DG(\tau)^+$, and
$\L_*$ is a positive endomorphism of $DG(\tau)$, which will be an order automorphism if, for
example, $\tau$ is surjective.  In this case, $DG(\tau)$ can be viewed as a $\Z[t,t^{-1}]$ module.

We now summarize this paper.  Sections 2--3 lead up to  the definition of the 
dimension group for a \pwm\ map
(Definition \ref{1.22}). Sections 4--6 develop basic properties
of the dimension group of a \pwm\ map, for example, characterizing when they are simple, and describing
a canonical set of generators for the dimension module (Theorem \ref{1.41}), which is quite useful in computing dimension
groups. Sections 7-11 compute the dimension group for various
families of interval maps, some of which are sketched below. Section 12 gives a dynamical description of a
C*-algebra $A_\tau$ such that $K_0(A_\tau)$ is isomorphic to $DG(\tau)$. 

\begin{figure}[htb]
\centerline{\includegraphics{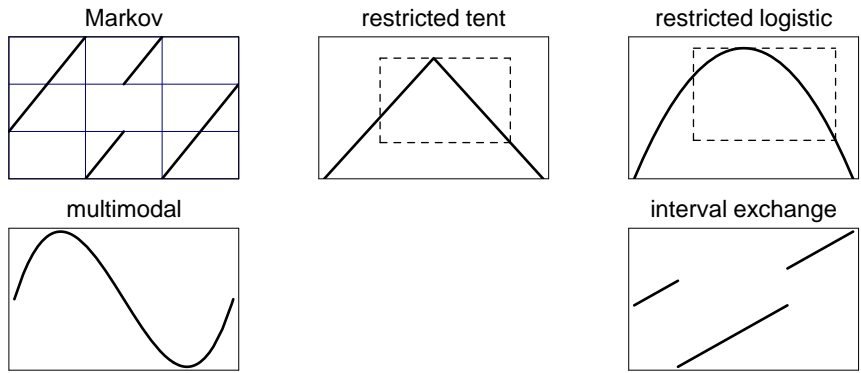}}
\caption{\label{fig1} Examples of piecewise monotonic maps}
\end{figure}

This paper initiates a program to make use of dimension groups (and more generally, C*-algebras and their K-theory) to find
invariants for interval maps.   
 In
\cite{ShuTrans},  we investigate dimension groups for transitive \pwm\ maps.  For such maps,  we 
describe the order on the dimension group in concrete terms, and show  that for some families of maps, the dimension triple
is a complete invariant for conjugacy.

In the  paper
\cite{DeaShu}, two C*-algebras $F_\tau$ and $O_\tau$ are associated with \pwm\ maps, and the properties of
these algebras are related to the dynamics. These algebras are analogs of the algebras $F_A$ and $O_A$
associated with shifts of finite type by Cuntz and Krieger \cite{CK}, and $F_\tau$ is isomorphic to the
algebra $A_\tau$ defined in Section 12 of this paper.

The author gratefully acknowledges partial support for this work by a Brachman Hoffman 
Fellowship from Wellesley College.

\section{Piecewise monotonic maps and local homeo\-morphisms}

 Let $I = [0,1]$. A map $\tau:I \to I$ is  {\it piecewise monotonic\/} if there are points $0 = a_0 < a_1 <
\ldots < a_n = 1$ such that $\tau|_{(a_{i-1},a_i)}$ is continuous and strictly monotonic for $1
\le i \le n$. We denote by
$\tau_i:[a_{i-1},a_i] \to I$ the unique continuous extension of $\tau|_{(a_{i-1},a_i)}$; note that
$\tau_i$ will be a homeo\-morphism onto its range.  We will refer to the ordered set $C = \{a_0, a_1, \ldots,
a_n\}$ as the {\it partition associated with $\tau$},  or as the endpoints of
the intervals of monotonicity. We will say this partition is {\it maximal\/} if the intervals
$(a_{i-1},a_i)$ are the largest open intervals on which
$\tau$ is continuous and strictly monotonic.  For each piecewise monotonic function there is a unique maximal
partition.  We do not assume a partition associated with $\tau$ is maximal unless that is specifically stated.
(We will be primarily interested in maximal partitions, but there will be a few cases where non-maximal
partitions are useful.)  

Our goal in this section is to modify $\tau$ slightly so that it becomes a local homeo\-morphism $\sigma$ on a larger space $X$, in such a
way that properties of $(I,\tau)$ and $(X,\sigma)$ are closely related. In general, points in the partition
$C$ cause trouble, since $\tau$ is usually neither locally injective nor an open map at such points, and may
be discontinuous. We will disconnect $I$ at points in $C$ (and at points in their forward and backward
orbit). 

We have to be a little careful about the
values of $\tau$ at points $c$ in $C$ where $\tau$ is discontinuous: what are relevant for our purposes are the left and right
limits $\lim_{x\to c^\pm} \tau(x)$, not the actual values $\tau(c)$.
If $\tau:I\to I$ is piecewise monotonic, we define a (possibly multivalued) function $\htau$ on $I$ by setting $\htau(x)$ to
be the set of left and right limits of $\tau$ at $x$.  At points where $\tau$ is continuous, $\htau(x) = \{\tau(x)\}$, and we
identify $\htau(x)$ with $\tau(x)$.  If
$A\subset I$, then
$\htau(A) =
\cup_{x\in A}
\htau(x)$, and
$\htau^{-1}(A) =
\{x
\in I
\mid
\htau(x)\cap A \not= \emptyset\}$.  Observe that $\htau$ is a closed map, i.e., if $A$ is a closed subset of $I$, then $\htau(A)$ is
closed, since $\htau(A) = \cup_i \tau_i(A)$. 

 The {\it
generalized orbit\/} of $C$ is the smallest subset $I_1$ of $I$ containing $C$ and closed under $\htau$ and $\htau^{-1}$.  This
is the same as the smallest subset of $I$ closed under $\tau_i$ and $\tau_i^{-1}$ for $1 \le i \le n$.  We define $I_0 =
I\setminus I_1$.   Next we describe the space resulting from disconnecting
$I$ at points in
$I_1$.  (The term ``disconnecting" we have borrowed from Spielberg \cite{Spi}.)

\defi{1.1}{Let $I = [0,1]$, and  let $I_0$, $I_1$ be as above. 
The {\it disconnection of $I$ at points in
$I_1$\/} is the totally ordered set $X$ which consists of  a copy of $I$ with the usual ordering,  but with each point
$x
\in I_1\setminus
\{0,1\}$ replaced by two points $x^- < x^+$. We equip $X$ with the order topology, and define the {\it collapse map\/}
$\pi:X\to I$ by
$\pi(x^\pm) = x$ for $x \in I_1$, and $\pi(x) = x$ for $x \in I_0$. We write  $X_1 = \pi^{-1}(I_1)$, and $X_0 = \pi^{-1}(I_0)
=X\setminus X_1$.}

If $x \in I_1$, then $x^-$ is the smallest point in $\pi^{-1}(x)$, and $x^+$ is the largest. For $x \in I_0$,
it will be convenient to define $x^- = x^+ = \widehat x$, where $\widehat x$ is the unique preimage under
$\pi$ of $x$.
For any pair $a, b \in
X$, we write
$[a,b]_X$ for the order interval
$\{x
\in X
\mid a
\le x
\le  b\}$. More generally, if $J$ is any interval in $\R$, then we write $J_X$ instead of $J\cap X$. 

\prop{1.2}{Let $X$ be the disconnection of $I$ at points in $I_1$, and $\pi:X\to I$ the collapse map. Then $X$ is homeomorphic
to a compact subset of
$\R$, and
\begin{enumerate}
\item $\pi$ is  continuous and order preserving.\label{1.2i}
\item $I_0$ is dense in $I$, and $X_0$ is dense in
$X$.\label{1.2ii}
\item $\pi|_{X_0}$ is a homeo\-morphism from $X_0$ onto $I_0$.\label{1.2iii}
\item $X$ has no isolated points.\label{1.2iv}
\item If $a, b \in I_1$, then $[a^+,b^-]_X$ is clopen in $X$, and every clopen subset of $X$ is a
finite disjoint union of such order intervals.\label{1.2v}
\end{enumerate}}

\prooff{Fix a listing of the
elements of $I_1$, say $I_1 = \{x_1, x_2, \ldots\}$.
 Define $\widehat I$ to be
$I$ with each point $x_k$ in $I_1\setminus \{0,1\}$ replaced by
a pair of points $x_k^- < x_k^+$ and a gap of length $2^{-k}$ inserted between these
points.  Then $\widehat I$ is a compact subset of $\R$, and the order topology on $\widehat I$ coincides with the topology
inherited from $\R$. Clearly $\widehat I$ is order isomorphic to the set $X$ in Definition \ref{1.1}.  (Hereafter we will identify
$X$ with $\widehat I$.)

(\ref{1.2i}) The inverse image of each open interval under $\pi$ is an open order interval in $X$,
so
$\pi$ is continuous.  It follows at once from the  definitions of $\pi$, and of the order on $X$, that  $\pi$ is
order preserving.

(\ref{1.2ii}) The complement of any finite subset of $I$ is dense in $I$. By the
Baire category theorem, since $I_1$ is countable, then $I_0$ is dense in $I$, and this implies density
of
$X_0$ in $X$.    

(\ref{1.2iii}) To prove that $\pi$ restricted to $X_0$ is a
homeo\-morphism, let $B$ be a closed subset of $X$. We will show $\pi(B
\cap X_0) = \pi(B) \cap I_0$. If $x \in \pi(B) \cap I_0$, then $x$ has a preimage in $B$,
and since $x \in I_0$, by definition of $X_0$, that preimage is also in $X_0$.  Thus  $ \pi(B) \cap
I_0 \subset \pi(B
\cap X_0)$, and the opposite containment is clear.  Hence $\pi(B\cap X_0)$ is closed in $I_0$, so
 $\pi|_{X_0}$ is a continuous, bijective, closed map from $X_0$ onto
$X_0$, and thus is a homeo\-morphism.

(\ref{1.2iv}) Since $I_0$ is infinite and dense in $I$, no point of $I_0$ is isolated in $I_0$, so by
(\ref{1.2iii}) no point of $X_0$ is isolated in
$X_0$. By density of $X_0$ in $X$, no point of $X$ is isolated.

(\ref{1.2v}) Observe that if $a \in I_1\setminus \{1\}$, then $a^+$ is not a limit from the left,
and if $b \in I_1\setminus \{0\}$,   then $b^-$ is not a limit from the right. Thus for $a, b \in
I_1$ with $a < b$, the
 set $[a^+,b^-]_X$ is clopen. Now let $E$ be any clopen subset of $X$, and let $x \in E$. Since
$E$ is clopen, for each point $x$ in $E$ we can find an open order interval  containing $x$ and 
contained in $E$.
The union of these intervals forms an
open cover of the compact set $E$, so there is a finite subcover. The union of overlapping open order intervals is again an
open order interval, so the intervals in this subcover can be expressed as a finite disjoint union of
open order intervals. By construction, not only do these cover $E$, but they are
contained in $E$, so $E$ is equal to their disjoint union. Since each of these
intervals is open, and each is the complement (in $E$)  of the union of the others,
each is also closed. Thus $E$ is a disjoint union of a finite collection of clopen order
intervals. 

We will be done if we show that for each clopen order interval $[c,d]_X$, there exist $a, b \in X_1$ with $c = a^+$ and $d = b^-$. By
density of $I_0$ in $I$, each point of $I_0$ is a limit from both the left and right of sequences from $I_0$,
  so each point of $X_0$ is a limit from both the left and right of sequences from $X_0$. Since
$[c,d]_X$ is open,
$c$ cannot be a limit from the left, and $d$ cannot be a limit from the right, 
so neither $c$ nor $d$ is in $X_0$. If $a\in I_1\setminus \{0\} $
and
$c = a^-$, then $c$ would be a left limit, which is impossible.  Thus $c = a^+$ for some $a \in I_1$, and
similarly $d = b^-$.}

  We now turn to constructing a local homeo\-morphism
$\sigma:X\to
X$.  Recall that if $a_0<a_1 < \cdots < a_n$ is the partition associated with a piecewise monotonic map $\tau$, we denote by
$\tau_i$ the extension of $\tau|_{(a_{i-1},a_i)}$ to a homeo\-morphism on $[a_{i-1},a_i]$.

\theo{1.3}{Let $\tau:I\to I$ be a piecewise
monotonic map,   $C = \{a_0, a_1, \ldots, a_n\}$ the endpoints of the
intervals of monotonicity, 
$I_1$ the generalized orbit of
$C$, and
$X$ the disconnection of $I$ at the points of $I_1$.   
\begin{enumerate}
\item The sets $\wJ_i = [a_{i-1}^+, a_i^-]_X$ for $i = 1, \ldots, n$ form a
partition of
$X$ into clopen sets.\label{1.3i} 
\item There is a unique continuous map
$\sigma:X\to X$ such that $\pi\circ \sigma = \tau \circ \pi$ on $X_0$.\label{1.3ii}
\item For  $1\le i \le n$, $\sigma$ is a  monotone homeo\-morphism from
$\wJ_i$ onto the clopen order interval with endpoints $\sigma(a_{i-1}^+)$ and
$\sigma(a_i^-)$.\label{1.3iii}
\item For $1\le i \le n$, $\pi\circ \sigma = \tau_i\circ \pi$ on $\wJ_i$.\label{1.3iv}
\end{enumerate}}

\prooff{For $1 \le i \le n$, let $J_i = [a_{i-1},a_i]$, and define $\wJ_i$ as in (\ref{1.3i}) above.
By Proposition \ref{1.2}, the order intervals
$\wJ_1, \ldots, \wJ_n$ are   clopen and  form a partition of $X$.  

Now we define $\sigma$ on $\wJ_1, \ldots,
\wJ_n$.   Fix $i$, and recall that
$\tau_i$ is strictly monotonic. We will assume that $\tau_i$ is increasing; the
decreasing case is similar. Then $\tau_i$  is a homeo\-morphism from
$J_i$ onto $[\tau_i(a_{i-1}), \tau_i(a_i)]$.

 Since
 $\pi(\wJ_i) = J_i$,  each point in
$\wJ_i$ has the form $x^+$ or $x^-$ for some $x \in J_i$.   If $y \in \wJ_i$ and $y = x^\pm$, then we define $\sigma(y) =
\sigma(x^\pm) = \tau_i(x)^\pm$. This is well defined, since if $y = x^+ = x^-$, then $x \in I_0$, so $\tau_i(x) \in I_0$, and
thus
$\tau_i(x)^+ = \tau_i(x)^-$.   It follows that $\sigma$ is an order
isomorphism from $\wJ_i$ onto $[\tau_i(a_{i-1})^+, \tau_i(a_{i})^-]_X$. Since the
topology on compact subsets of $\R$ coincides with the order topology, $\sigma$ is an
increasing homeo\-morphism from $\wJ_i$ onto $\sigma(\wJ_i)$, so (\ref{1.3iii}) holds. Since
$\tau_i(I_1)
\subset I_1$, then by Proposition
\ref{1.2},
$\sigma(\wJ_i) = [\tau_i(a_{i-1})^+,
\tau_i(a_i)^-]_X$ is clopen. 
 Since each order interval $\wJ_i$ is clopen, then $\sigma$
is continuous on all of $X$. From the definition of $\sigma$, $\pi\circ \sigma = \tau \circ \pi$ on $X_0$, and
so by continuity of $\sigma$ and density of $X_0$, (\ref{1.3ii}) follows. Now (\ref{1.3iv}) is
immediate.}

From the proof above, if $\tau$ is increasing on the interval  $J_i = [a_{i-1},a_i]$,  and if $c, d \in J_i\cap I_1$ with $c < d$,
then
\begin{equation}
\sigma([c^+,d^-]_X) = [\tau_i(c)^+,\tau_i(d)^-]_X,\label{(1.1)}
\end{equation}
while if $\tau$ decreases on $J_i$ then
\begin{equation}
\sigma([c^+,d^-]_X) = [\tau_i(d)^+,\tau_i(c)^-]_X.\label{(1.2)}
\end{equation}

\cor{1.4}{$\pi|_{X_0}$ is a conjugacy from
$(X_0,\sigma)$ to $(I_0,\tau)$.}

\prooff{By Proposition \ref{1.2}, $\pi|_{X_0}$ is a homeo\-morphism onto $I_0$, and by
Theorem \ref{1.3}, $\pi$ intertwines 
$(X_0,\sigma)$ and $(I_0,\tau)$.}

The map $\sigma$ in Theorem \ref{1.3} is a local homeo\-morphism, i.e., it is an open map, and each point admits an open neighborhood
on which $\sigma$ is a homeo\-morphism. 

\defi{1.5}{We will call the map $\sigma:X\to X$ in Theorem \ref{1.3} the  {\it local homeo\-morphism associated
with
$\tau$}.}

\exe{1.6}{Let $\tau:[0,1] \to [0,1]$ be the full tent map  given by $\tau(x)
=1-|2x-1|$. Here the set $C$ of endpoints of intervals of monotonicity is $\{0, 1/2, 1\}$. This set is invariant under $\tau$, so the
generalized orbit $I_1$ is 
$\cup_{n\ge0}
\tau^{-n}(C)$, which is just the set of dyadic rationals. The order intervals $\wJ_1, \wJ_2$ defined in Theorem \ref{1.3} satisfy
$\sigma(\wJ_1) =
\sigma(\wJ_2) = X$. For $x \in X$, let $S(x) = s_0s_1s_2 \ldots$, where $s_n = i$ if $\sigma^nx \in \wJ_i$. Then $S$ is a
conjugacy from $(X,\sigma)$ onto the full (one-sided) 2-shift.}

For a general \pwm\ map $\tau:I\to I$, the associated local homeo\-morphism $\sigma:X\to X$ will not be a shift of finite type, nor even a
subshift.  For example, for a logistic map
$\tau(x) = kx(1-x)$ with an attractive fixed point,  $X$ will contain non-trivial connected components. 
\smallskip

In the remainder of this section, we will see that properties of a \pwm\ map $\tau:I\to I$ and the associated local homeo\-morphism $\sigma$
are closely related. We illustrate this relationship for transitivity.

\defi{1.7}{If $X$ is any topological space, and $f:X\to X$ is
a continuous map, then
$f$ is {\it transitive\/} if for each pair
$U, V$ of non-empty open sets, there exists $n \ge 0$ such that
$f^n(U)
\cap V
\not=
\emptyset$. We say $f$ is {\it strongly transitive\/} if for
every non-empty open set $U$, there exists $n$ such that
$\cup_{k=0}^nf^k(U) = X$.}

\defi{1.7.1}{
 If $\tau:I\to I$ is \pwm, we view $\tau$ as undefined at the set $C$ 
of endpoints of intervals of monotonicity, and say $\tau$ is transitive
if  for each pair
$U, V$ of non-empty open sets, there exists $n \ge 0$ such that
$\tau^n(U)
\cap V
\not=
\emptyset$. We say $\tau$ is {\it strongly transitive\/} if for
every non-empty open set $U$, there exists $n$ such that $\cup_{k=0}^n\htau^k(U) = I$.
(Recall that
$\htau$ denotes the (possibly multivalued) function whose value
at each point $x$ is given by 
 the left and right hand limits of $\tau$ at $x$.)}

 It is well known that for $X$ compact metric with no isolated
points, and $f:X\to X$ continuous, transitivity is equivalent to the existence of a point with a dense orbit. The same
equivalence holds if $\tau:I\to I$ is \pwm, and is viewed as undefined on the set $C$ of  endpoints of intervals of monotonicity.
(Apply the standard Baire category argument to show that transitivity implies the existence of a dense orbit, cf. e.g. 
\cite[Thm. 5.9]{WalBook}).) 

If $\tau:I\to I$ is \pwm\ and continuous, then any dense orbit must eventually stay out of $C$. It follows that 
 the definitions of transitivity in Definitions \ref{1.7} and \ref{1.7.1} are consistent.  In addition, in this case $\tau = \htau$,
so the definitions of strong transitivity also are consistent.

For general compact metric spaces $X$, transitive maps need not be strongly
transitive.  For example, the (two-sided) shift on the space of all bi-infinite sequences of
two symbols is transitive in the sense of Definition \ref{1.7}, but is not strongly transitive, since
the complement of a fixed point is open and invariant. However, the next proposition shows
that every piecewise monotonic transitive map is strongly transitive. (For continuous \pwm\ maps, this follows from \cite[Thm.
2.5]{Preston}.) In the proof below, we will repeatedly use the fact that the collapse map $\pi:X\to I$, when restricted to
$X_0$, is a conjugacy from $(X_0,\sigma)$ onto $(I_0,\tau)$, cf. Corollary \ref{1.4}. 

\prop{1.8}{For a \pwm\ map $\tau:I\to I$, with associated local homeo\-morphism $\sigma:X\to X$,
the following are equivalent.
\begin{enumerate}
\item $\tau$ is transitive.
\item $\tau$ is strongly transitive.
\item $\sigma$ is transitive.
\item $\sigma$ is strongly transitive.
\end{enumerate}}

\prooff{Assume $\tau$ is transitive. View $\tau$ as undefined on the set of endpoints
of intervals of monotonicity. Let
$V$ be an open interval. Choose an open interval
$W$ such that
$\overline{W} \subset V$. By \cite[Prop. 2.6]{BuzziSpec},  or \cite[Corollary on p. 382]{HofPWI}, there exists an $N$ such that 
$\cup_0^\infty \tau^k(W) =  \cup_0^N \tau^k(W)$. The right side is a finite union of open
intervals. Since $\tau$ is transitive, the left side is dense,
and thus  $I\setminus \cup_0^N \tau^k(W)$ is finite.  On the
other hand, for each $k$, $\tau^k(W) 
 \subset \htau^k(\overline{W}) \subset \htau^k(V)$. 
Since $\htau$ is a closed map, it follows that
$\cup_0^N\htau^k(\overline{W}) = I$, so $\cup_0^N\htau^k(V) =
I$. Thus $\tau$ is strongly transitive. 

Strong transitivity of $\tau$ is equivalent to strong
transitivity of $\tau|_{I_0}$.  This in turn is equivalent to
strong transitivity of  $\sigma|_{X_0}$, and then to strong
transitivity of $\sigma$ on $X$.  Finally, if $\sigma$ is
transitive on $X$, it is transitive on $X_0$, so $\tau$ is
transitive on $I_0$, which in turn implies transitivity of
$\tau$ on $I$.}

\section{Dimension groups} 

In this section we will associate a dimension group $DG(\tau)$ with each piecewise monotonic map $\tau$. Our method will be
to define a dimension group $G_\sigma$, where $\sigma$ is the local homeo\-morphism associated with $\tau$, and then define
$DG(\tau) = G_\sigma$.   We start by defining a class of functions that includes  the local homeo\-morphisms associated with
piecewise monotonic maps, but also includes any local homeo\-morphism of a zero dimensional metric space.

\defi{1.10}{Let $X$ be a topological space. A map $\sigma: X\to X$ is a {\it piecewise
homeo\-morphism\/} if $\sigma$ is continuous and open, and $X$ admits a finite partition into clopen sets $X_1, X_2,
\ldots, X_n$ such that $\sigma$ is a homeo\-morphism from $X_i$ onto $\sigma(X_i)$ for $i = 1, \ldots, n$.}

If $\tau:I \to I$ is piecewise monotonic, and $\sigma:X\to X$ is the associated local homeo\-morphism, then $\sigma$ is
a piecewise homeo\-morphism, cf. Theorem \ref{1.3}.  
Every local homeo\-morphism
$\sigma$ on a compact zero dimensional metric space $X$ is a piecewise homeo\-morphism.  In fact, each open neighborhood of a
point in a zero dimensional space contains a clopen neighborhood of that point, so there is a cover of $X$ by clopen subsets
on which
$\sigma$ is injective.  Since $X$ is compact, we can take a finite subcover, and the resulting partition that arises
from all intersections of these subsets and their complements is the desired partition. Note that a composition of piecewise 
homeo\-morphisms is again a piecewise homeo\-morphism.

 Hofbauer \cite{HofPWI} and Keller \cite{Kel}  have studied the dynamics of maps that are {\it piecewise invertible}, i.e., that
meet the requirements of a piecewise homeo\-morphism except for the requirement that the map be open. 

We begin by defining $G_\sigma$ as an ordered abelian group; later we will show $G_\sigma$ is in fact a dimension group. Recall
that an abelian group $G$ is {\it ordered\/} if there is a subset  $G^+$ such that $G = G^+ - G^+$, $G^+ \cap
(-G^+) =
\{0\}$ and $G^+ + G^+ \subset G^+$. Then for $x, y  \in G$ with
 $y-x \in G^+$, we write $x \le y$. We write $n x$ for the sum of $n$ copies of $x$.  An element $u$ of
an ordered abelian group $G$ is an {\it order unit\/} if for each $x$ in $G$ there is a positive integer $n$
such that $-n u \le x \le n u$. 

The motivating idea for $G_\sigma$ is to build a group out of equivalence classes of clopen subsets of $X$. As discussed in
the introduction, this idea originated with Krieger \cite{KrDim}, and was extended by Boyle-Fiebig-Fiebig \cite{BFF}.  In
\cite{BFF} a dimension group is associated with any surjective local homeo\-morphism on a zero dimensional compact metric space.
Our context is a slight generalization of that in \cite{BFF}, namely, piecewise homeo\-morphisms on an arbitrary compact metric
space.

In what follows, a key role will be played by the {\it transfer map\/}.  We will  use
this map both to define the dimension group for a piecewise homeo\-morphism, and to provide a key tool in
computing this group.

\defi{1.11}{Let $X, Y$ be compact metric spaces. If $\sigma:Y\to X$ is any map which is finitely-many-to-one,
the {\it transfer map
$\L_\sigma$ associated with
$\sigma$\/} is defined by
\begin{equation}
(\L_\sigma f)(x) = \sum_{\sigma(y) =x} f(y)
\end{equation}
for each $f:X\to\C$. (We will usually have $Y = X$.)
We will write $\L$ instead of $\L_\sigma$ when the map $\sigma$ intended is clear from
the context. }

We note some simple properties of $\L$, which are immediate consequences of the definition. If $f\ge 0$, then $\L
f
\ge 0$.  If $f\ge 0$ and $\L f = 0$, then $f = 0$.
If 
$\sigma$ is injective on
$E$, then 
\begin{equation}
\L\chi_E =
\chi_{\sigma(E)}.\label{(1.3)}
\end{equation}
  If $\sigma:X\to X$ is a piecewise homeo\-morphism, then every function in $C(X,\Z)$ can
be written as a sum $\sum n_i\chi_{E_i}$ where $\sigma$ is injective on $E_i$, so 
$\L$ maps $C(X,\Z)$ into $C(X,\Z)$.
Observe that for $\sigma:X\to X$,
\begin{equation}
\L_{\sigma^n} = (\L_\sigma)^n.\label{(1.4)}
\end{equation}
If $X$ is a compact metric space, we will write $1_X$ (or simply $1$) for the function constantly 1 on $X$. 

\defi{1.12}{Let $X$ be a compact metric space, and $\sigma:X\to X$ a piecewise homeo\-morphism. Write $\L$ for
$\L_\sigma$. We define an equivalence relation ~ on $ C(X,\Z)$ by $f \sim g$ if $\L^n f = \L^n g$ for some $n \ge
0$, and write $G_\sigma$ for the set of equivalence classes, and $G_\sigma^+ = \{[f] \mid f \ge 0\}$. We
define addition on $G_\sigma$ by
$[f]+[g] = [f+g]$, and we order $G_\sigma$ by   $[f] \le [g]$ if $[g-f] \in G_\sigma^+$.  We call $[1_X]$ the {\it
distinguished order unit\/} of $G_\sigma$.}

It is easily verified that $G_\sigma$ is an ordered abelian group. 
We will show $G_\sigma$ is a dimension group, after developing some properties of the ordering on $G_\sigma$.
When $X$ is a zero dimensional space and $\sigma:X\to X$ is a local homeo\-morphism, then $G_\sigma$ is the same as the
``images group" defined by Boyle-Fiebig-Fiebig in \cite{BFF}.  Our definition is different from theirs, but equivalent, as observed in
 \cite[Rmk. 9.4]{BFF}. 

\lem{1.13}{Let $X$ be a compact metric space and $\sigma:X\to X$ a piecewise homeo\-morphism.  Then for each clopen
subset $F$ of $\sigma(X)$, there exists a clopen subset $E$ of $X$ such that $\sigma$ is injective on $E$ and
$\sigma(E) = F$.}

\prooff{Since $\sigma$ is a piecewise homeo\-morphism, there exists  a partition $ G_1, G_2, \ldots, G_n$ of $X$  into
clopen sets on which
$\sigma$ is injective.  If
$F
\subset \sigma(G_i)$ for some $i$, then  the set $E = G_i \cap \sigma^{-1}(F)$ has the desired properties. 
Otherwise, let
$F_1, F_2, \ldots, F_k$ be the partition of $F$ generated by the sets $\{F\cap \sigma(G_i)\mid 1 \le i \le n\}$. 
For each
$k$, choose a clopen set $E_k$ such that $\sigma$ is injective on $E_k$ and $\sigma(E_k) = F_k$. Then $\sigma$ is
injective on $E = E_1\cup E_2\cup \cdots E_k$, and $\sigma(E) = F$.}

If
$f
\in C(X,\Z)$, we define the {\it support\/} of
$f$ to be the set of $x$ in $X$ where $f$ is nonzero, and denote this set by $\supp
f$. Note that by the definition of $\L$ we have
\begin{equation}
\supp \L f \subset \sigma(\supp f),\label{(1.5)}
\end{equation}
with equality if $f \ge 0$.
In particular, since $\L_\sigma^n = \L_{\sigma^n}$, we have
\begin{equation}
\supp \L^nf \subset \sigma^n(X).\label{(1.6)}
\end{equation}
A result similar to the following can be found in \cite[Lemma 2.8]{BFF}.

\lem{1.14}{Let $X$ be a compact metric space, and
$\sigma:X\to X$ a piecewise homeo\-morphism. If $f \in C(X,\Z)$, and the support of $f$ is
contained in $\sigma^n(X)$, then there exists $f_0 \in C(X,Z)$
such that $\L^nf_0 = f$. If $f \ge 0$, then $f_0$ can be
chosen so that $f_0 \ge 0$. }

\prooff{Let $f \in C(X,\Z)$, with $\supp f \subset \sigma^n(X)$. Write $f = \sum_{i=1}^p n_i \chi_{E_i}$, with $E_1, \ldots,
E_p$ disjoint, and
$n_1,
\ldots, n_p
\in
\Z$. We may assume that for each index $i$, $n_i
\not= 0$. Then each $E_i$ is contained in the support of $f$, and
thus is contained in $\sigma^n(X)$. Applying Lemma \ref{1.13} to
$\sigma^n$, for each index $i$ we can find a clopen set $F_i$
such that
$\sigma^n$ is injective on $F_i$ and $\sigma^n(F_i) = E_i$. Let $f_0
=\sum n_i \chi_{F_i}$. Then
$\L^n f_0 = \sum n_i \chi_{\sigma^n(F_i)}= \sum n_i \chi_{E_i}
= f$. If $f \ge 0$, we can arrange $n_i> 0$ for all $i$, so $f_0 \ge 0$.}

\lem{1.15}{Let $X$ be a compact metric space, and
$\sigma:X\to X$ a piecewise homeo\-morphism. In $G_\sigma$, $[f] \le [g]$ iff $\L^nf \le \L^n g$
for some $n \ge 0$.}

\prooff{It suffices to show $[f] \ge [0]$ iff $\L^nf \ge 0$
for some $n \ge 0$.  If $[f] \ge [0]$, then by definition there is $h \ge 0$ such that $f\sim h$, and therefore
$\L^n f = \L^n h \ge 0$ for some $n\ge 0$. Conversely, fix $n$ and suppose  $\L^nf \ge 0$. By (\ref{(1.6)}) and Lemma \ref{1.14}, there
exists
$h
\ge 0$ such that $\L^n f = \L^n h$. Then $f\sim h$, so $[f] \ge [0]$.}

\defi{1.16}{$\L_*:G_\sigma
\to G_\sigma$ is defined by $\L_*[f] = [\L f]$.}

This is an injective homomorphism. If $\sigma$ is surjective, by Lemma \ref{1.14} $\L$ is surjective, so $\L_*$ is surjective, and
thus is an automorphism of $G_\sigma$.   Note, however, that $\L_*$ is usually not unital.

If $\sigma$ is injective on a clopen set $E$, then by (\ref{(1.3)})
\begin{equation}
\L_*[\chi_E] = [\chi_{\sigma(E)}].\label{(1.7)}
\end{equation}

\notate{Let $G_1 \mapright{T_1} G_2 \mapright{T_2}G_3 \cdots$ be a sequence of ordered abelian groups with positive connecting
homomorphisms.  We review the construction of the inductive limit of this sequence. Let $G_\infty$ be the set of sequences
$(g_1, g_2,
\ldots)$ such that
$g_i \in G_i$ for each $i$, and such that there exists $N$ such that $T_ng_n =g_{n+1}$ for all $n \ge N$.  Make
$G_\infty$ into a group with coordinate wise addition. Consider two sequences $(g_i)$, $(h_i)$ to be equivalent if they
eventually agree, i.e., if there exists $N$ such that $g_i = h_i$ for $i \ge N$.  The quotient group $G$ is then an
abelian group. Denote equivalence classes by square brackets.  Define
$G^+$ to be the set of equivalence classes of sequences  that are eventually positive, and order $G$ by $g \le h$ if $h-g \in
G^+$. Then
$G$ is the inductive limit in the category of ordered abelian groups.   If $G_1$, $G_2$, $\ldots$ are identical, and $T_1$,
$T_2$, $\ldots$ are the same map,  we call $G$ the {\it stationary\/} inductive limit with  connecting maps $T_1:G_1 \to G_1$.}

\lem{1.17}{Let $X$ be a compact metric space, and
$\sigma:X\to X$ a piecewise homeo\-morphism.  Then $G_\sigma$ is isomorphic as an ordered group to the inductive limit
\begin{equation}
C(X,\Z)\mapright{\L}\L(C(X,\Z)\mapright{\L}\L^2C(X,\Z)\mapright{\L}\cdots.
\end{equation}
and this isomorphism carries the map $\L_*$ to the shift map $[(g_1, g_2, \ldots)] \mapsto [(g_2, g_3, \ldots)]$.}

\prooff{Let $G$ be the inductive limit of this sequence. Every element of $G$ has the
form
$[(0_0,0_1,\ldots,0_{n-1}, \L^ng, \L^{n+1}g, \ldots)] = [(g, \L g, \ldots)]$, where $0_i$ denotes a zero in the $i$-th
position. The map $\pi:G\to G_\sigma$ given by
$\pi([(g,\L g, \ldots)]) = [g]$ is an order isomorphism of $G$ onto $G_\sigma$, and carries $\L_*$ to the shift map.}

If $\sigma:X\to X$ is surjective, then $\L:C(X,\Z) \to C(X,\Z)$ is also surjective (Lemma \ref{1.14}). Hence for surjective
$\sigma$, the dimension group $G_\sigma$ is isomorphic to the inductive limit of the sequence
\begin{equation}
C(X,\Z)\mapright{\L}C(X,\Z)\mapright{\L}C(X,\Z)\mapright{\L}\cdots.\label{(1.8)}
\end{equation}

\lem{1.18}{Let $X$ be a compact metric space, and $\sigma:X\to X$ a piecewise homeo\-morphism. For each $n \ge 0$,
$\L^nC(X,\Z)$ (with the order inherited from $C(X,\Z)$)  is isomorphic as an ordered group to
$C(\sigma^n(X),\Z)$.}

\prooff{By (\ref{(1.6)}), $\supp \L^n f \subset \sigma^n(X)$, so the map $f \mapsto
f|_{\sigma^n(X)}$ is an order isomorphism from $\L^n C(X,\Z)$ into $C(\sigma^n(X),\Z)$. To prove
this map is surjective, let $g \in C(\sigma^n(X),\Z)$. Extend $g$ to be zero on
$X\setminus \sigma^n(X)$, so that $g$ is in $C(X,\Z)$ with support in $\sigma^n(X)$. There exists $f \in C(X,\Z)$ with
$L^n f = g$ by Lemma \ref{1.14}. Then $\L^n f |_{\sigma^n X} = g$, so $f \mapsto f|_{\sigma^n X}$ is surjective.}

\defi{1.19}{An ordered abelian group $G$ has the {\it Riesz
decomposition property\/} if whenever $x_1, x_2, y_1, y_2 \in G$ with $x_i \le y_j$ for
$i,j = 1,2$, then there exists $z \in G$ with $x_i \le z \le y_j$ for
$i,j = 1,2$.  $G$ is {\it unperforated\/} if $ng \ge 0$ for some $n > 0$ implies $g \ge
0$. A {\it dimension group\/} is a countable ordered abelian group
$G$ which is unperforated and which has  the Riesz decomposition property.}

A more common definition is that a 
dimension group is the inductive limit of a sequence of ordered groups $Z^{n_i}$
(where 
$Z^{n_i}$ is given the usual coordinate-wise order.)  By a result of Effros,
Handelman, and Shen \cite{EHS}, this is equivalent to the definition we have given.  Two standard references on dimension groups
are the books of Goodearl \cite{Go} and of Effros \cite{Eff}.

If $G$ is a dimension group with
a distinguished order unit $u$, then we refer to the pair $(G,u)$ as a {\it unital dimension group\/}. A
homomorphism between unital dimension groups $(G_1, u_1)$, $(G_2,u_2)$ is {\it unital\/} if it takes $u_1$ to
$u_2$. 

Below $C(X,\Z)$ is
viewed as an ordered abelian group with the usual  addition and ordering of functions.

\lem{1.20}{If $X$ is a compact metric
space, then $C(X,\Z)$ is a dimension group, with order unit $1_X$.}

\prooff{We show $G_\sigma$ is countable.
Since $X$ is a
compact metric space, then the space $C(X)$ of continuous complex valued functions
on $X$ is a separable Banach space.  Any two characteristic functions of clopen
sets are a distance $2$ apart in $C(X)$, so there are at most
countably many clopen sets in $X$. Since every element of $C(X,\Z)$ is a finite
integral combination of characteristic subsets of clopen sets, then $C(X,\Z)$ is also
countable.

Suppose that $f_i, g_j$ are in $C(X,\Z)$ with $f_i \le g_j$ for $i=
1, 2$,
$j = 1, 2$.  If
$h =
\max(f_1,f_2)$, then
$f_i\le h
\le g_j$ for $i= 1, 2$, $j = 1, 2$, so $C(X,\Z)$ satisfies the Riesz decomposition property. It is evident that $C(X,\Z)$ is
unperforated.}

\cor{1.21}{If  $X$ is a compact metric space, and
$\sigma:X\to X$ a piecewise homeo\-morphism, then $G_\sigma$ is a dimension group, with distinguished order unit $[1_X]$.}

\prooff{Since $C(X,\Z)$ is countable (cf. Lemma \ref{1.20}), and since $G_\sigma$ is a quotient of $C(X,\Z)$, then $G_\sigma$ is
countable. With the order inherited from
$C(X,\Z)$, 
$\L^nC(X,\Z)$ is a dimension group, since it is isomorphic to the dimension group $C(\sigma^n(X),\Z)$, cf. Lemma \ref{1.18}. Each
of the properties defining a dimension group (being unperforated, and the Riesz decomposition property) are preserved by
inductive limits, so by Lemma \ref{1.17}, $G_\sigma$ is a dimension group. For any $f \in C(X,\Z)$ there exists $k$ such that $-k
1_X \le f
\le k 1_X$, so
$[1_X]$ is an order unit.}

\defi{1.22}{If $\tau:I \to I$ is piecewise monotonic, and $\sigma:X\to X$ is the associated local homeo\-morphism,
then 
$DG(\tau)$ is defined to be the dimension group $G_\sigma$. Here we assume the partition $C$ associated with $\tau$ is the maximal one;
if not, then we instead write $DG(\tau, C)$ for this dimension group.}

\section{Reduction to surjective maps}

In this section, we will see that we can often reduce the computation of the dimension group $DG(\tau)$ to the case where $\tau$ is
surjective.

\defi{1.23}{If $X$ is a  metric space and $\sigma:X\to X$ is continuous, we say $\sigma$ is {\it eventually
surjective} if
$\sigma^{n+1}(X) =
\sigma^n(X)$ for some $n \ge 0$. In this case we refer to $\sigma^n(X)$ as the {\it
eventual range\/} of $\sigma$. If $\tau:I\to I$ is piecewise monotonic but not continuous, we say $\tau$ is eventually surjective
if $\tau|_{I_0}$ is eventually surjective, and  $\tau$ is surjective if $\htau$ is surjective (i.e., if $\htau(I) = I$).}

If $\tau$ is \pwm\ and continuous, it is straightforward to check that the two definitions of eventual surjectivity in Definition
\ref{1.23} that are applicable are consistent.

\lem{1.23.1}{If $\tau:I\to I$ is \pwm, and $\sigma:X\to X$ is the associated local homeo\-morphism, then $\sigma$
is essentially surjective iff $\tau$ is essentially surjective, and $\sigma$ is surjective iff $\tau$ is surjective.}

\prooff{Since $I_0$ is dense in $I$, and $\htau$ is a closed map, surjectivity of $\tau$ is equivalent to surjectivity of
$\tau|{I_0}$. By density of
$X_0$ in
$X$ and continuity of
$\sigma$, surjectivity of $\sigma$ is equivalent to surjectivity of $\sigma|_{X_0}$. By the
conjugacy of $(I_0,\tau)$ and
$(X_0,\sigma)$, surjectivity of $\tau$ is equivalent to surjectivity of $\sigma$. 

Similarly,  $\tau^{n+1}(I_0) =
\tau^n(I_0)$ is equivalent to $\sigma^{n+1}(X_0) =
\sigma^n(X_0)$, which in turn is equivalent to $\sigma^{n+1}X = \sigma^n X$, so eventual surjectivity of $\sigma$ is equivalent to
eventual surjectivity of $\tau$.}

An alternative definition for eventual surjectivity of discontinuous $\tau$ might be 
that
$\htau^{n+1}(I) =\htau^n(I)$ for some $n \ge 0$. This implies
$\tau^{n+1}(I_0) =
\tau^n(I_0)$, so
$\tau$ is eventually surjective in the sense of Definition \ref{1.23}. The converse is also true, but is a little tedious to prove when
$\tau$ is discontinuous, and won't be needed in the sequel, so we have chosen to define eventual surjectivity of $\tau$ as in
Definition
\ref{1.23} instead.

\prop{1.24}{If $X$ is a compact metric space and $\sigma:X\to X$ is an eventually surjective piecewise homeo\-morphism with
eventual range $Y$, then $G_\sigma$ is isomorphic to the stationary inductive limit
\begin{equation}
C(Y,\Z) \mapright{\L}
C(Y,\Z)\mapright{\L}
C(Y,\Z)\mapright{\L}\cdots.
\end{equation}}

\prooff{This follows from Lemmas \ref{1.17} and \ref{1.18}.}

\exe{1.25}{If $\tau:[0,1]
\to [0,1]$ is defined by $\tau(x) = kx(1-x)$, and $2 \le k <
4$, then $\tau([0,1]) = \tau([0,1/2]) = [0,k/4] \supset
[0,1/2]$, so $\tau^2([0,1]) = \tau([0,1])$. Thus $\tau$ is
eventually surjective, but is not surjective. By Lemma \ref{1.23.1}, the associated piecewise homeo\-morphism $\sigma:X\to X$ is also
eventually surjective but not surjective. On the other hand, if $0 \le k <
2$, then the sets $\tau^n([0,1])$ are nested with intersection
$\{0\}$, so $\tau$ is not eventually surjective.}

The following is \cite[Thm. 4.3]{BFF}, adapted to our current context.

\prop{1.26}{Let $X$ be a compact metric space, and
$\sigma:X\to X$ a piecewise homeo\-morphism. Then $\L_*$
is a homomorphism and an order isomorphism of $G_\sigma$ onto its image, and is surjective iff $\sigma$ is eventually
surjective.}

\prooff{We observed after Definition \ref{1.16} that $\L_*$ is injective. Since $\L$ is a positive operator, it follows that
$\L_*$ is positive. If $\L_*[f]\ge 0$, then $[\L f]
\ge 0$, which by definition implies $\L^{n+1}f\ge 0$ for some $n \ge 0$, and thus $[f] \ge 0$. Hence
$\L_*$ is an order isomorphism onto its image.

Next we show $\L_*$ is surjective if $\sigma$ is eventually surjective.
Suppose that $\sigma^{n+1}(X) = \sigma^n(X)$, and let $f \in
C[X,\Z]$. Then the support of
$\L^n f$ is contained in $\sigma^n(X) = \sigma^{n+1}(X)$, so there
exists $g\in C(X,\Z)$ such that $\L^{n+1}g = \L^n f$. It
follows that $\L_* [g] = [f]$.

Finally, suppose $\L_*$ is surjective. Choose $f \in C(X,\Z)$ such that $\L_*[f] = [1]$. Then for some $n \ge 0$, $\L^{n+1}f = L^n 1$.
Then $\sigma^n X = \supp \L^n 1 = \supp \L^{n+1}f \subset \sigma^{n+1}X$, so $\sigma^nX = \sigma^{n+1}X$. Thus $\sigma$ is eventually
surjective.}

\defi{1.27}{Let $X$ be a compact metric space, and
$\sigma:X\to X$ a piecewise homeo\-morphism. Then $(G_\sigma, G_\sigma^+, (\L_\sigma)_*)$ is the {\it dimension triple\/}
associated with $\sigma$.}

Suppose $X_1, X_2$ are compact metric spaces,  $\sigma_i:X_i\to X_i$ for $i = 1, 2$ are piecewise  homeo\-morphisms, and
$\psi:X_1\to X_2$ is a piecewise homeo\-morphism that intertwines $\sigma_1$ and $\sigma_2$, i.e., $\psi\circ \sigma_1 =
\sigma_2\circ
\psi$. Then $\L_\psi\L_{\sigma_1} =
\L_{\psi\circ \sigma_1} = \L_{\sigma_2\circ \psi} = \L_{\sigma_2}\L_\psi$.  It follows that $\L_\psi\L_{\sigma_1}^n =
\L_{\sigma_2}^n\L_\psi$. Thus the map
$(\L_\psi)_*: G_{\sigma_1} \to G_{\sigma_2}$ defined by $(\L_\psi)_*[f] = [\L_\psi f]$ is well defined, and is a positive
homomorphism that intertwines $(\L_{\sigma_1})_*$ and $(\L_{\sigma_2})_*$.

\prop{1.28}{Let $X$ be a compact metric space, and
$\sigma:X\to X$ an eventually surjective piecewise homeo\-morphism. Let $Y$ be the
eventual range of $\sigma$, and choose $N$ so that $\sigma^N (X) = Y$. Then there is a group and order isomorphism $\Phi$
from
$G_\sigma$ onto $G_{\sigma|_Y}$, which carries the automorphism $(\L_\sigma)_*$
to the automorphism $(\L_{\sigma|_Y})_*$, and carries the distinguished order unit $[1_X]$ to  the order unit $[(\L^N
1)|_Y]$.}

\prooff{Let $\psi:Y\to X$ be the inclusion map; note that $\psi$ intertwines $\sigma|_Y$ and $\sigma$. By the remarks preceding
this proposition,
$\Phi_0 =
(\L_\psi)_*$ is a positive homomorphism from $G_{\sigma|_Y}$ into $G_\sigma$, intertwining $(\L_{\sigma|_Y})_*$ and
$(\L_\sigma)_*$.  We will show that
$\Phi_0$ is a  order isomorphism of $G_{\sigma|_Y}$ onto $G_\sigma$.

Observe that $\L_\psi f$ is the function that agrees with $f$ on $Y$, and is zero off $Y$; in particular, $\L_\psi$ is 1-1 on
$C(Y,\Z)$, and $\L_\psi f \ge 0$ iff $f \ge 0$.  

To prove that $\Phi_0$ is surjective, let $f \in C(X,\Z)$. Note that $\supp \L_\sigma^N f \subset
\sigma^N(X) = Y$. Since $\sigma$ is surjective on $Y$,  Lemma \ref{1.14} (applied to $\sigma|_Y$) implies that there exists $g \in
C(Y,\Z)$ such that $\L_{\sigma|_Y}^N g = (\L_\sigma^N f)|_Y$. Applying
$\L_\psi$ to both sides gives $\L_\sigma^N \L_\psi g = \L_\psi(\L_\sigma^Nf|_Y)=
\L_\sigma^N f$  (since $\supp \L_\sigma^Nf \subset Y$), so $\Phi_0([g]) = [\L_\psi g] = [f]$. Thus $\Phi_0$ is surjective. Furthermore,
if $[f]= [\L_\psi g] \ge 0$,  then there exists $n$ such that $\L_\sigma^n \L_\psi g \ge 0$. Then $\L_\psi L_{\sigma|_Y}^n g \ge 0$,
which implies
$L_{\sigma|_Y}^n g \ge 0$. Thus $[g] \ge 0$, so $\Phi_0$ is an order isomorphism.

Finally, if $f = 1_X$ in the last paragraph,  and $g\in C(Y,\Z)$ is chosen so that  $\L_{\sigma|_Y}^N g = (\L_\sigma^N
1_X)|_Y$, then $\Phi_0([g]) = [1_X]$. Thus $\Phi_0^{-1}[1_X] =
(\L_{\sigma|_Y})_*^{-N}[(\L_\sigma^N 1_X)|_Y]$. Since $\sigma|_Y$ is surjective, by Proposition
\ref{1.26},   $(\L_{\sigma|_Y}^N)_*$ is an order automorphism of $C(Y,\Z)$. Now the map $\Phi  =
(\L_{(\sigma|_Y)})_*^N\Phi_0^{-1}$ satisfies the  conclusion of the proposition.}

\noindent{\bf Remark}
Suppose that $\tau:I\to I$ is continuous, and piecewise  monotonic, and let $C$ be the set of endpoints of the associated
maximal partition. If $\tau$ is eventually surjective, with eventual range $J$, then $\tau|_J$ is surjective, and we can reduce
calculation of
$DG(\tau)$ to finding $DG(\tau|_J,C')$ for a suitable partition $C'$ as follows.

Choose $N$ so that $J = \tau^N(X) = \tau^{N+1}(X)$. Since $\tau$ is continuous, then $J$ is a
closed interval.  Let
$C'=(C\cap J)\cup \tau^N(C)$. (If desired, we may omit from $C'$ any
points in $\tau^N (C)$ that are  in the generalized orbit of $C\cap J$ under $\tau|_J$.) 
Associate with
$\tau|_J$ the partition whose endpoints are the points in
$C'$. Then the generalized orbit of
$C'$ under $\tau|_J$ will be $I_1 \cap J$. If $\sigma:X\to X$ is the local homeo\-morphism associated with
$\tau$, then $\sigma$ will be eventually surjective; let $Y$ be its eventual range.  Then the local
homeo\-morphism associated with $\tau|_J$ (with the partition $C'$) will be $\sigma|_Y$. Thus by Proposition
\ref{1.28}, $ DG(\tau)  \cong DG(\tau|_J, C')$.

\exe{1.29}{Let $\tau$ be the map  on $[0,2]$ which 
has the values $\tau(0) = 0$, $\tau(1/2) = 1$, $\tau(1) = 0$, $\tau(3/2) = k < 1$, $\tau(2) = 0$, with $\tau$ linear
between these points.  The eventual range of $\tau$ is $J = \tau([0,2]) = [0,1]$.  We take as a  partition of
$\tau|_J$ the non-maximal partition $C'=\{0 , 1/2 , k, 1\}$.  Then $DG(\tau) \cong DG(\tau|_J, C')$.}

\exe{1.29.1}{Let $\tau$ be the logistic map in Example \ref{1.25}, with $2 \le k < 4$. The eventual range of $\tau$ is $J =
\tau([0,1]) = [0,\tau 1]$. We can take as partition for $\tau|_J$ the maximal partition $C' = \{0,1/2, \tau 1\}$. Thus $DG(\tau)
\cong DG(\tau|J)$.}

\section{Simplicity of the dimension group}

An {\it order ideal\/} in a dimension group $G$ is a subgroup $J$ such that $a, b \in J$
and
$a
\le c \le b$ implies $c \in J$. An ideal $J$ of a dimension group is {\it positively
generated\/} if $J = J^+ - J^+$ where $J^+ = J\cap G^+$. A dimension group $G$ is
{\it simple\/} if it has no  positively generated order ideals other than $\{0\}$ and $G$. We will use the
following equivalent condition:
$G$ is simple iff every nonzero positive element of $G$ is an order unit.
We now examine conditions that guarantee that $G_\sigma$ is simple.

\defi{1.30}{Let $X$ be a compact metric space and $\sigma:X\to X$ a continuous map.
Then $\sigma$ is {\it (topologically) exact\/} if for every open set $V$ there exists $n \ge 0$ such
that
$\sigma^n(V) = X$. If $\tau:I\to I$ is piecewise monotonic (but not necessarily continuous everywhere), then
$\tau$ is topologically exact if for  every open set $V$ there exists $n \ge 0$ such
that
$\htau^n(V) = I$.  (Recall that if $\tau$ is continuous, then $\htau = \tau$, so these two definitions of exactness
are consistent.)}

If $X$ is a compact metric space, and $\sigma:X\to X$ is topologically exact, then evidently $\sigma$ is topologically
mixing. The converse is not necessarily true, but is  true for a broad class of \pwm\ maps.  For a more
thorough discussion of the relationship, see \cite[Prop. 4.9]{ShuTrans}. 

Two examples of topologically exact maps are $\beta$-transformations, and tent maps with
slopes
$\pm s$ with
$\sqrt{2} < s \le 2$, restricted to a suitable interval \cite{Jonk-Rand}.

\lem{1.31}{Let $\tau:I\to I$ be piecewise  monotonic, with associated local homeo\-morphism $\sigma:X\to X$. Then $\tau$
is topologically exact iff $\sigma$ is topologically exact.}

\prooff{We first show that $\tau$ is exact iff $\tau|_{I_0}$ is exact. Assume $\tau$ is exact, and let $V$ be an open subset
of $I_0$. Choose an open subset $W$ of $I$ such that $W\cap I_0 = V$. By exactness of $\tau$, there exists $n$ such that
$\htau^n(W) = I$. Since $I_0$ is forward and backward invariant under $\htau$, then $\tau^n(V) = \tau^n(W\cap I_0) =
I_0$, which proves that $\tau|_{I_0}$ is exact.

Conversely, suppose $\tau|_{I_0}$ is exact, and let $W$ be a non-empty open subset of $I$. Choose a non-empty open set $V
\subset W$ such that $\overline{V} \subset W$. 
Then $V\cap I_0$ is open in $I_0$, and is
non-empty by density of $I_0$. Choose $n$ so that $\tau^n(V\cap I_0) = I_0$. Then $\htau^n(V) \supset
\tau^n(V\cap I_0) = I_0$. Since $\htau$ is a closed map, then $\htau^n(\overline{V})$ is a closed subset of $I$,
containing the dense set $I_0$, so equals $I$. Then $\htau^n(W) = I$, which proves that $\htau$ is exact.}

Recall that a compact metric space is zero dimensional if there is a base of clopen sets, and that this is equivalent to
being totally disconnected, cf. \cite{HW}.

\theo{1.32}{Let $X$ be a compact metric space, and
$\sigma:X\to X$ a piecewise homeo\-morphism. Then $G_\sigma$ is simple
iff 
\begin{enumerate}
\item $\sigma$ is eventually surjective, and\label{1.32i}
\item for every compact open subset $V$ of the eventual range $Y$ of $\sigma$,
there   exists
$n \ge 0$ such that $\sigma^n(V) = Y$. \label{1.32ii}
\end{enumerate}
\noindent In particular, if $\sigma$ is eventually
surjective and is topologically exact when restricted to its eventual range, then
$G_\sigma$ is simple. The converse is true if $X$ is totally disconnected.}

\prooff{Suppose that (\ref{1.32i}) and (\ref{1.32ii}) hold.   By Proposition \ref{1.28}, if $Y$ is
the eventual range of $\sigma$, then $G_\sigma$ is isomorphic to
$G_{\sigma|_Y}$, so we may assume without loss of generality that $\sigma$ is surjective.  
Let
$f
\in C(X,\Z)^+$, with $f \not= 0$, and let $E$ be a clopen set such that $\chi_E\le f$. (Since $f\ge 0$ takes on integer values, such an
$E$ exists.) Choose
$n
\ge 0$ such that
$\sigma^n(E) = X$, and choose a positive integer $k$ such that $\L^n 1 \le k 1$. (Since $\L^n 1$
is bounded, such an integer $k$ exists.)  Then
\begin{equation}
\L_\sigma^nf \ge \L_{\sigma^n}( \chi_E) \ge \chi_{\sigma^n(E)} =  1_X \ge (1/k)\L^n 1_X.
\end{equation}
Thus $k[f] \ge  [1_X]$. Since $[1_X]$ is an order unit, then so is $[f]$, and this proves that $G_\sigma$ is simple.

Conversely, suppose $G_\sigma$ is simple. We first show that $\sigma$ is eventually
surjective. We have
$X
\supset
\sigma(X)
\supset
\sigma^2(X)
\cdots$. Since $G_\sigma$ is simple, then there is some $k \ge 0$ such that $[1_X] \le
k[\chi_{\sigma(X)}]$.  Then for some $n \ge 0$, $\L^n 1_X \le k \L^n
\chi_{\sigma(X)}$. The support of $\L^n 1_X$ is $\sigma^n(X)$, and the
support of $\L^n \chi_{\sigma(X)}$ is $\sigma^{n+1}(X)$, so $\sigma^n(X) \subset
\sigma^{n+1}(X)$. Since the opposite inclusion also holds, we must have equality, so
$\sigma$ is eventually surjective.

Now let $Y$ be the eventual range of $\sigma$. If $Y \not= X$, replace $\sigma$ by $\sigma|_Y$, so that we may
assume $\sigma$ is surjective. Let $E$ be any non-empty clopen subset of
$X$. By simplicity of $G_\sigma$, there exists $k > 0$ such that $k[\chi_E] \ge [1_X]$, and thus  there exists $n \ge 0$ such that
$k \L^n\chi_E
\ge  \L^n 1_X$. Comparing supports, we conclude that $\sigma^n(E) = X$, so (\ref{1.32ii}) holds.

Finally, suppose that $X$ is totally disconnected, and that $G_\sigma$ is simple. Then by
(\ref{1.32i}), $\sigma$ is eventually surjective. Let $Y$ be the eventual range of  $\sigma$.  Let
$V$ be any non-empty open subset of $Y$. Then $V$ contains a non-empty clopen subset $W$, and so by
(\ref{1.32ii}) there exists
$n$ such that $\sigma^n(W) = Y$. Therefore $\sigma^n(V) = Y$, proving that $\sigma$ is topologically exact on $Y$.}

\cor{1.33}{If $\tau:I\to I$ is piecewise monotonic and topologically exact, then $DG(\tau)$ is a simple
dimension group.}

\prooff{Let $\sigma:X\to X$ be the associated local  homeo\-morphism, and assume $\tau$ is topologically exact.  By Lemma
\ref{1.31}, $\sigma$ is topologically exact, so by Theorem \ref{1.32}, $DG(\tau)$ is simple.}

 By Theorem \ref{1.32}, simplicity of $DG(\tau)$ is most informative when $X$ is totally disconnected. We now
explore when this occurs. Recall that the generalized orbit of a subset $A$ of $I$ is the smallest set containing $A$ and closed
under
$\htau$ and
$\htau^{-1}$. 

\lem{1.34}{Let $\tau:I\to I$ be piecewise  monotonic,  let $C$ be the set of points in the associated partition, and let
$\sigma:X\to X$ be the  associated local homeo\-morphism.  Then
$X$ is totally disconnected iff the generalized orbit $I_1$ of $C$ is dense in $I$. }

\prooff{Suppose that $I_1$ is dense
in
$I$. Let $d \in X$, and let $(a,b)_X$ contain $d$. We will assume $\pi(a) < \pi(d) < \pi(b)$; the cases where $\pi(a) =
\pi(d)$ or $\pi(d) = \pi(b)$ can be proved in a similar manner. By density of $I_1$, we can choose $x,  y \in I_1$ with $\pi(a) < x
< \pi(d) < y < \pi(b)$. Then $a < x^+ < d < y^- < b$, and $[x^+, y^-]_X$ is a clopen order interval containing $d$ and
contained in $(a,b)_X$. Thus the clopen sets form a base for the topology of $X$, so $X$ is totally disconnected.

Conversely, suppose that $I_1$ is not dense in $I$. Let $(a,b)$ be an open interval in $I$ containing
no point of $I_1$. Then every point in $(a,b)$ is in $I_0= I\setminus I_1$, so every point in  the open set
$\pi^{-1}(a,b)$ is in $X_0$. Since each clopen subset of $X$ contains points of
$X_1$ (cf. Proposition \ref{1.2}), then 
$\pi^{-1}(a,b)$  is a non-empty open subset of
$X$ containing no non-empty clopen subset.  Thus $X$ is not totally disconnected.}

 Below ``interval" will always
denote an interval that is not a single point. 

\defi{1.36}{Let $\tau:I\to I$ be \pwm, and
$J$ an open interval. Then $J$ is a {\it homterval\/} if  $\tau^n$ is a homeo\-morphism on $J$ for all $n$.}

Note that if $J$ is a
homterval for $\tau$, then for each
$p \ge 0$, $J$ will be a homterval for $\tau^p$. For additional background on homtervals, see the books of Collet-Eckmann
\cite[\S II.5]{ColEck} or de Melo-van Strien \cite[Lemma II.3.1]{Melo-S}.  The following result characterizes the existence of
homtervals for polynomials.

Let $\tau:I\to I$ be piecewise monotonic and continuous. 
An {\it attracting periodic orbit\/} is
a periodic orbit such that the set of points whose orbits converge to that
periodic orbit has non-empty interior.

\prop{1.36.1}{If $\tau:I\to I$ is a polynomial of degree $\ge 2$, then $\tau$ has a homterval iff $\tau$ has an attractive periodic
orbit.}

\proof{If  $f$ has a homterval, by \cite[Lemma II.3.1 and Thm.~A, p.267]{Melo-S}, there exists
an attractive periodic orbit.
Conversely, let $x_0$ be an attractive periodic point of $f$, say $f^p(x_0) = x_0$.  Let
$h = f^{2p}$; then $h$  has $x_0$ an attractive fixed point. Note $h'(x_0) =
((f^p)'(x_0))^2
\ge 0$. If $h'(x_0) > 0$,   since
$x_0$ is attractive, 
there must be an open interval
$V$ with endpoint
$x_0$ such that $0 < h'< 1$ on $V$. By the mean value theorem, $h(V) \subset V$, and $h$ is 1-1 on $V$. Thus
$V$ is a homterval for $h$, and therefore also for $f$.

 If $h'(x_0) = 0$, choose $\delta > 0$ such that $0 < |h'(x)| < 1$ for $0 < |x-x_0| < \delta$. Let $W = (x_0-\delta, x_0 +
\delta)\cap I$. By the mean value theorem, $h$ maps $W\setminus \{x_0\}$ into itself.  If $V$ is a component of
$W\setminus
\{x_0\}$, then
$V$ is a homterval for $h$, and thus for $f$.}

\prop{1.37}{Let $\tau:I\to I$ be \pwm, and let $\sigma:X\to X$ be the associated local
homeo\-morphism. The following are equivalent.
\begin{enumerate}
\item $I$ contains no homterval.\label{item1.37(i)}
\item $I_1$ is dense in $I$.\label{item1.37(ii)}
\item $X$ is totally disconnected. \label{item1.37(iii)}
\end{enumerate}
In particular, these hold if $\tau$ is transitive.}

\prooff{(\ref{item1.37(i)}) $\Rightarrow$ (\ref{item1.37(ii)}) If $I$ contains no homterval, then for any open interval $J$,
there exists
$n$ such that
$\tau^n(J)$ meets $C$. Thus $I_1$ is dense in $I$.

(\ref{item1.37(ii)}) $\Leftrightarrow$ (\ref{item1.37(iii)}) Lemma \ref{1.34}.

(\ref{item1.37(ii)}) $\Rightarrow$ (\ref{item1.37(i)})  Suppose  that there exists a homterval $J$; we will show $I_1$ is not
dense in
$I$. We first show that either (1) there exists $p \ge 1$ and a
 homterval whose images under $\tau^p$ are disjoint, or (2) there exists $p \ge 1$ and an interval fixed
pointwise by $\tau^p$. 

Suppose the images of $J$ under $\tau$ are not disjoint, so that
$\tau^n(J)
\cap
\tau^{n+p}(J)
\not=
\emptyset$ for some $n \ge 0$ and $p \ge 1$. Then for $k \ge 0$, $\tau^{n+kp}(J) \cap \tau^{n+(k+1)p}(J) \not=
\emptyset$, so $V = \cup_{k=0}^\infty \tau^{n+kp}(J)$ is an interval invariant
under $\tau^p$. Since $J$ is a $\tau^p$-homterval, therefore $\tau^p$ is continuous and monotonic on each of the overlapping
open intervals
$\tau^{n+kp}(J)$ for $k = 0, 1, \ldots$, so $\tau^p$ is a homeo\-morphism from $V$ into $V$. Replacing $ \tau^p$ by $\tau^{2p}$
if necessary, we may assume that $\tau^p$ is increasing on $V$.

If some point $x$ in $V$ is
not fixed by $\tau^p$,  then the $\tau^p$-orbit of $x$ is a monotonic sequence
converging to some point of $\overline{V}$.  Let $K$ be the open interval with endpoints $x$ and $\tau^p(x)$.
Then
$K$ is an $\tau^p$-homterval whose  images under $\tau^p$ are disjoint. Thus we've  shown either (1) or (2) must hold. We will 
show in both cases that $I_1$ is not dense in $I$. 

Let $p\ge 1$ and $h = \htau^p$, and suppose that $J$ is a
homterval whose images under $h$ are disjoint. We will prove that $h^N(J) \cap I_1$ is finite for some $N \ge 0$, showing that $I_1$ is
not dense in $I$.  Let $D = \cup_{k=0}^{p}\htau^k(C)$.  Observe that
$\cup_{k=1}^\infty h^k(D) = \cup_{n=1}^\infty \htau^n(C)$ is invariant under $\htau$, and that the generalized orbit of $C$
under $\tau$ is the same as the generalized orbit of $D$ under $h$, i.e., 
\begin{equation}\label{eqn16}
I_1 = \{x \mid \htau^mx \cap \htau^n C \not= \emptyset \hbox{ for some $m, n \ge 0$}\} =  \{x \mid  h^jx
\cap h^k D
\not= \emptyset \hbox{ for some $j, k \ge 0$}\}.
\end{equation} In the rest of this proof,
``orbit" means the orbit under the map $h$ unless otherwise specified.

Since the orbit of $J$ consists of disjoint intervals, by replacing $J$ by $h^k(J)$ for suitable $k \ge 0$,
we can arrange that the orbit of $J$ is disjoint from the finite set $C\cup h(C)$.  (Since $\htau$ may be multivalued at points in $C$,
it simplifies the discussion to avoid such points where possible.)

Let $D_1$ consist of those points in $D$ whose orbits miss $C$ and meet the orbit of $J$. If $d \in D_1$, let $n_d$ be the
least integer such that $\hbox{orbit}(d) \cap \tau^{n_d}(J) \not=
\emptyset$, and let $N = \max_{d\in D_1} n_d$, and $J_N = \tau^N(J)$. (If $D_1 = \emptyset$, define $J_N = J$.) The orbit of
each point
$d$ in
$D_1$ meets
$J_{n_d}$, and therefore meets $J_N$. Furthermore, $\hbox{orbit}(d)$ meets $J_N$ in a unique point $y_d$, since having $\tau^kd
\in J_N$ and
$\tau^{k+q}d\in J_N$ would imply $\tau^qJ_N\cap J_N \not= \emptyset$, contrary to the assumption that the orbit of $J$ consists
of disjoint intervals.

Let $P = J_N \cap I_1$; we will prove $P$ is finite. If $x \in P$, by (\ref{eqn16}) the orbit of $x$ meets the orbit
of  some
$d_0
\in D$. Then there is a sequence $d_0,
d_1,
\ldots, d_k$ such that
$d_{i+1}\in h(d_i)$ for
$0\le i\le k-1$, and
$d_k
\in
\orbit(x)$. By construction of $J$, the orbit of $J$ misses $C$, so the orbit of $x$ misses $C$; thus $d_k \notin C$. Let $j$ be
the least index such that $d_j, d_{j+1}, \ldots, d_k$ miss $C$, and set $d = d_j$. Then $d \in D_1$, and the orbit of
$d$ meets the orbit of
$x$.
 Let $y_d$ be the unique point where the orbit of
$d$ meets
$J_N$. Distinct points in
$J_N$ have disjoint orbits, so $y_d = x$. Thus for each $x \in P$, there exists a point $d \in D_1$ whose orbit meets $J_N$ exactly in
the single point $x$. Since $D_1$ is finite, then $P$ is finite, which finishes the proof that $I_1$ is not dense in $I$.

Suppose instead that there exists an interval $V$ fixed pointwise by $h = \tau^p$ for some $p \ge 1$.  Replacing $V$
by a component of $V\setminus (C\cup \htau C)$ if necessary, we can assume that $V$ is disjoint from $C\cup \htau C$. Let
$D_2$ be the set of points in $D$ whose orbits land in $V$ and miss $C$. Each orbit of a point in $D_2$ is finite; let $A$ be
the union of such orbits.  Then no point in $C$ has an orbit that meets $V\setminus A$, so any component of $V\setminus A$ is
disjoint from $I_1$.
 Thus $I_1$ is not dense in
$I$.

Finally, suppose $\tau$ is transitive. Then $\tau$ is strongly transitive by Proposition \ref{1.8}. If $V$ is any non-empty open
subset of $I$, and $c \in C$, by strong transitivity there exists $n\ge 0$ such that $c \in\htau^n(V)$. Then $V$ meets
$\htau^{-n}(c)$, so $V\cap I_1 \not= \emptyset$. Thus $I_1$ is dense in $I$.}

\exe{1.35}{Let $1 < s \le 2$ and let $\tau$ be the tent map with slopes $\pm s$.  The length of each interval not
containing the critical point 1/2 is expanded by $\tau$ by a factor $s$, so for every open interval $J$, there exists $n$
such that 1/2 is in $\tau^n(J)$. Thus $I_1$ is dense in $I$, so if $\sigma:X\to X$ is the local homeo\-morphism associated
with $\tau$, then $X$ is totally disconnected.}

We have the following partial converse for Corollary \ref{1.33}. 

\cor{1.38}{If $\tau:I\to I$ is piecewise monotonic,
is surjective, has no homtervals, and $DG(\tau)$ is simple, then $\tau$ is topologically exact.}

\prooff{Let $\sigma:X\to X$ be the associated local homeo\-morphism. Then $X$ is totally disconnected 
(Proposition \ref{1.37}), and by hypothesis, the dimension
group
$G_\sigma= DG(\tau)$ is simple. Since $\tau$ is surjective, then $\sigma$ is surjective, cf. Lemma \ref{1.23.1}. By Theorem
\ref{1.32}, $\sigma$ is topologically exact on its eventual range, namely, $X$, and then by Lemma \ref{1.31},
$\tau$ is topologically exact.}

\exe{1.39}{Let $\tau:I\to I$ be a polynomial of degree $\ge 2$, and let $\sigma:X\to X$ be the associated local homeo\-morphism. Then
$X$ will be totally disconnected iff $\tau$ has no attractive periodic orbit (Propositions \ref{1.36.1} and \ref{1.37}).}

\section{Module structure of the dimension group}

In this section we will generally work with eventually surjective piecewise monotonic maps $\tau:I\to I$,  
cf. Definition \ref{1.23}. For such maps, if $\sigma:X\to X$ is the associated local homeo\-morphism, then by
Proposition \ref{1.26}, 
$(\L_\sigma)_*$ is an order automorphism of $DG(\tau)= G_\sigma$. We use this to give $DG(\tau)$
the structure of a module. We will show that
this module is finitely generated, and will identify generators.

\defi{1.40}{Let $\tau:I\to I$ be  piecewise
monotonic and eventually surjective.   Then we view
$DG(\tau)$ as a
$\Z[t,t^{-1}]$-module by defining 
\begin{equation}
\left(\sum_{i=-n}^m z_it^i\right)[f] = \left(\sum_{i=-n}^m z_i(\L_*)^i\right)[f]
\end{equation}
for $z_{-n}, \ldots, z_m \in \Z$ and $[f] \in DG(\tau)$.  Thus $DG(\tau)$ is both a dimension group and a
$\Z[t,t^{-1}]$ module, and we call it the {\it dimension module\/} for $\tau$.}

Let $0 = a_0 < a_1 < \ldots < a_n = 1$ be the partition associated with $\tau$. Recall that $\tau_i$ denotes the unique
extension of
$\tau|_{(a_{i-1},a_i)}$ to a homeo\-morphism from $[a_{i-1},a_i]$ onto 
its image. 

If $a, b \in I_1$ with $a < b$, then we write $I(a,b)$ for the order interval
$[a^+,b^-]_X$. If $a > b$, then we define $I(a,b) = I(b,a)$, and we set $I(a,a) = 0$. Recall that each such set $I(a,b)$
is clopen, and that every clopen subset of
$X$ is a finite disjoint union of such sets, cf. Proposition \ref{1.2}.  We generally will not distinguish between an
order interval $E$ in $X$, its characteristic function $\chi_E \in C(X,\Z)$, and the equivalence class
$[\chi_E] \in
DG(\tau)$; it should be clear from the context which is intended. In particular, if $a, b \in I_1$, we consider $I(a,b)$ as a
member of $DG(\tau)$.

We also observe that if $a, b$ are in  $[a_{i-1},a_i]\cap I_1$,
then  by equations (\ref{(1.1)}) and (\ref{(1.2)}),
\begin{equation}
\sigma( I(a,b)) = I(\tau_i(a),\tau_i(b)).     \label{(1.9)}
\end{equation}

If $D$ is a finite subset of $\R$, we will say distinct points $x, y \in D$ are
{\it adjacent\/} in $D$ if there is no element of $D$ between them. 

\theo{1.41}{Let $\tau:I\to I$ be \pwm\ and eventually surjective, with associated partition $C = \{a_0, a_1,
\ldots, a_n\}$.  Let $M$ be any member of the set $\htau C = \{\tau_1(a_0),
\tau_1(a_1), \tau_2(a_1),\tau_2(a_2)\ldots, 
\tau_n(a_{n-1}), \tau_n(a_n)\}$.
Let
\begin{equation}
\J_1 =
\{\,I(c,d)
\mid \hbox{$c$, $d$ are adjacent points in 
$\{\,a_0, a_1,
\ldots, a_n, M\,\}$}\,\},
\end{equation}
and let $\J_2$ be the set of intervals corresponding to jumps at partition points, i.e., 
\begin{equation}
\J_2 = \{\,I(\tau_i(a_i), \tau_{i+1}(a_i)) \mid 1 \le i
\le n-1\,\}.
\end{equation}
 Then
$DG(\tau)$ is generated as a module by $\J_1 \cup \J_2$.}

\prooff{Let
$H$ be the submodule of
$DG(\tau)$ generated by $\J_1\cup \J_2$. By Proposition \ref{1.2},
$DG(\tau)$ is generated by the intervals $I(c,d)$ for $c,d \in
I_1$, so it suffices to show each such interval is in $H$. Define $c
\sim d$ if
$I(c,d) \in H$. Then we need to show
$c\sim d$ for all $c,d \in I_1$.  This relation is clearly reflexive and
symmetric. To show it is transitive, suppose
$c \sim d$ and $d
\sim e$. If $c \le d \le e$, then $I(c,d)$ and $I(d,e)$ are in $H$, and are disjoint intervals, so the  sum
$I(c,e)$ of their corresponding characteristic functions is also in
$H$, and thus
$c \sim e$. If $c\le e \le d$, then $I(c,e) = I(c,d) - I(e,d) \in H$, so again $c \sim e$. The case $e \le c\le
d$ is similar.

By the definition of $\J_1$, each pair of adjacent points in $\{\,a_0, a_1,
\ldots, a_n, M\,\}$ is equivalent. By transitivity, all points in $\{\,a_0, a_1,
\ldots, a_n, M\,\}$ are equivalent. 
Suppose
 $c,d \in [a_{i-1},a_i]\cap I_1$. We will show 
\begin{equation}
c\sim d \iff \tau_i(c) \sim
\tau_i(d).\label{(1.10)}
\end{equation}
Recall from Theorem \ref{1.3} that $\sigma$ is a homeo\-morphism on
$[a_{i-1}^+,a_i^-]_X
\supset I(c,d)$, so by (\ref{(1.7)}) we have
$\L_*(I(c,d)) =
\sigma(I(c,d))$.  From (\ref{(1.9)}),  $I(\tau_i(c),\tau_i(d))=\L_*(I(c,d))$.  Thus 
$I(\tau_i(c),\tau_i(d))\in H$ iff $I(c,d) \in H$, which gives (\ref{(1.10)}).

In particular, $\tau_i(a_{i-1}) \sim \tau_i(a_i)$ for $1 \le i \le n$. By the
definition of $\J_2$, we have $\tau_i(a_i)\sim\tau_{i+1}(a_i)$ for $1 \le i \le n-1$. Then by
transitivity we conclude that $\tau_1(a_0),
\tau_1(a_1), 
\ldots,  \tau_n(a_{n-1}), \tau_n(a_n)$ are equivalent. By definition, $M$ is one
of these numbers, so we conclude that each of these is also equivalent to each of
$a_0, a_1, \ldots,a_n$. 

Let $P$ be the set of all points in $I_1$
equivalent to $a_0$. So far we have shown that $P$ contains $a_0, a_1, \ldots,a_n$, and 
$\tau_1(a_0),
\tau_1(a_1),
\ldots, \tau_n(a_{n-1}), \tau_n(a_n)$. Let $x \in P$.  Then for some index $i$,
$x$ is in $[a_{i-1},a_i]$. Since $x \in P$, then $x \sim a_0$, and since $a_0 \sim
a_{i-1}$, then $x \sim a_{i-1}$. By (\ref{(1.10)}), $\tau_i(x) \sim \tau_i(a_{i-1}) \sim
a_0$, so
$\tau_i(x)
\in P$. Thus $P$ is closed under application of each $\tau_i$.

Now suppose that $x
\in P$ and $y\in [a_{i-1},a_i]$  with $\tau_i(y)
= x$. Since $\tau_i(a_i) \in P$, then $\tau_i(y) = x \sim \tau_i(a_i)$. By (\ref{(1.10)}), $y \sim a_i$, so $y \in P$.
 Therefore $P$ contains
$\{a_0, a_1, \ldots, a_n\}$ and is closed under application of  each $\tau_i$ and
$\tau_i^{-1}$, so must coincide with $I_1$. }

\cor{1.42}{If  $\tau:I\to I$ is a
continuous, surjective piecewise monotonic map with associated partition
$\{a_0, a_1, \ldots, a_n\}$, then
$DG(\tau)$ is generated as a module by $I(a_0,a_1)$, $I(a_1,a_2)$, $\ldots$,
$I(a_{n-1},a_n)$.}

\prooff{Since $\tau$ is continuous, then all intervals in $\J_2$ in Theorem
\ref{1.41} are empty, so can be omitted. Since $\tau$ is continuous and monotonic
on each interval $[a_i,a_{i+1}]$, its maximum must occur at one of $a_0,
\ldots, a_n$. On the other hand, since $\tau$ is surjective, its maximum
value must be $a_n$. Thus $a_n = \tau(a_j)$ for some $j$, so $M$ as defined in
Theorem \ref{1.41} is equal to $a_n$. Now the result follows from Theorem \ref{1.41}.}

\exe{1.43}{Let $1 \le s \le 2$, and let $\tau:[0,1]\to [0,1]$ be the map
$\tau(x) = s x$ for
$0
\le x
\le 1/2$ and $\tau(x) = s - s x$ for $x > 1/2$.  By Theorem \ref{1.41}, for the   maximal
partition   $\{0,1/2, 1\}$, taking $M = \tau(0) = 0$, $DG(\tau)$ is
generated as a $Z[t,t^{-1}]$ module by $I(0,1/2)$ and $I(1/2, 1)$. Furthermore, we have
$\L_*(I(0,1/2)) = I(0,s/2) = \L_*I(1/2,1)$, so in fact this module is  singly generated,
with generator $I(0,s/2)$.}

\exe{1.44}{Let $\tau$ be the map $\tau(x) = 2x$ for $0 \le x \le 1/3$ and $\tau(x)
= (4/3 - x)$ for $1/3 < x \le 1$. Then $DG(\tau)$ is generated by the intervals
$I(0,1/3)$, $I(1/3, 1)$, and $I(2/3,1)$. The last interval is included because of the jump
discontinuity at $1/3$; note that $\tau_1(1/3) = 2/3$ while $\tau_2(1/3) = 1$.}

\section{Cyclic dimension modules}

We will see later that the dimension modules of surjective unimodal maps are cyclic (i.e., singly generated). In the current
section, we find conditions that guarantee that cyclic dimension modules are free. These conditions  involve requirements that
certain orbits of critical points are infinite and disjoint.\label{Cyclic}

Let $\tau:I\to I$ be \pwm, with associated local homeo\-morphism $\sigma:X\to X$. 
We will
identify functions in $C(X,\Z)$ with functions in $L^1(\R)$ in the following manner. If $g \in C(X,\Z)$, there is a
unique
$f\in L^1(\R)$ such that $(f\circ \pi)(x) = g(x)$ for $x \in X_1$, with $f(x) = 0$ for $x \notin I$.
  We will always choose $f$ to be continuous except
at a finite set of points in $I_1$. (This is possible by Proposition \ref{1.2}(\ref{1.2v}). For
example, for
$a, b \in I_1$, the function in $L^1(\R)$ corresponding to the characteristic function of the order
interval $[a^+,b^-]_X$ is the characteristic function of $[a,b]$.)  Any two such choices for $f$
will agree except on a finite set.  

Since we identify functions that agree except at a finite number of
points, then by a ``discontinuity" for a function in $C(X,\Z)$, viewed as a function on $\R$, we will mean a
point where left and right limits exist but are different. We denote the set of discontinuities of $g$ by $\D(g)$. Note that
$\D(1_X) = \{0,1\}$.

Next we will develop some basic facts about $\D(g)$. For the definition of $I(a,b)$, see the beginning of
the previous section.

 Let $\tau:I\to I$ be \pwm, with associated partition $C$, and let $C_0 = C\setminus \{0,1\}$. For $g \in
C(X,\Z)$, let $x_0  < x_1 < \cdots < x_q $ be the points in $C_0$, together with the points
where $g$ is discontinuous.  Then for some $n_1, \ldots, n_q$ in $\Z$, we can write 
\begin{equation}
g = \sum_{i=1}^q n_i I(x_{i-1},x_i),\label{(1.11)}
\end{equation}
and
\begin{equation}
\L g= \sum_i n_i\sigma(I(x_{i-1},x_i))=\sum_i n_iI(\tau_i(x_{i-1}),\tau_i(x_i)). \label{(1.12)}
\end{equation}
It follows that
\begin{equation}
\D(\L g) \subset \htau(\D(g)) \cup \htau C_0,
\end{equation}
so for all $n \ge 0$,
\begin{equation}
\D(\L^n g)\subset \htau^n\D(g) \cup \bigcup_{k=1}^n\htau^k C_0.\label{(1.13)}
\end{equation}

Note in Lemma \ref{1.45}, we interpret $\tau e \notin \htau C_0$, where $C_0 = C\setminus\{0,1\}$, to mean in
particular that
$e
\notin C_0$, so that
$\tau$ is single-valued at $e$.

\lem{1.45}{Let $\tau:I\to I$ be \pwm,  with associated partition $C$, and let $C_0 = C\setminus \{0,1\}$.  Let
$e\in I_1$ satisfy $\tau e\notin \htau C_0$.  If $g \in C(X,\Z)$ is discontinuous at $e$,  and
$g$ is continuous at all points in $\tau^{-1}(\tau e)$ except $e$, then $\L g$ is discontinuous
at
$\tau e$.}

\prooff{Since we are viewing functions in $C(X,\Z)$ as members of $L^1(\R)$, if $a, b \in I_1$,  we will write $\tau([a,b])$ for
characteristic function of the closed interval $\overline{\tau((a,b))}$. With the notation in (\ref{(1.11)}), since
$\D(g)
\subset
\{x_0,
\ldots, x_q\}$,  there is some index
$j$ such that
$e = x_j$. First consider the case
$0 < j < q$.   Since
$g$ is discontinuous at $e$, then  $n_j \not= n_{j+1}$. Since $e \notin C_0$, then $\tau$ is
monotonic and continuous on $(x_{j-1},x_{j+1})$. Therefore $\tau([x_{j-1},x_j])$ and
$\tau([x_j,x_{j+1}]$ are closed intervals with just the endpoint $\tau x_j = \tau e$ in common.  Since
$n_j\not= n_{j+1}$,  then $
n_j\tau([x_{j-1},x_j])+n_{j+1}\tau([x_j,x_{j+1}])$ is discontinuous at
$\tau x_j = \tau e$. On the other hand, for all $i \not= j$, either $x_i \in C_0$ or else $x_i$ is a
point of discontinuity of
$g$. If $x_i \in C_0$, then $ \tau e \notin \htau C_0$ implies that $\tau e \notin \htau x_i$. If $x_i \in \D(g)$, by hypothesis $x_i
\notin \tau^{-1}\tau e$, so again $\tau e \notin \htau x_i$. 
Thus none of the endpoints of the intervals
$\tau([x_{i-1},x_i])$ for $i\notin \{j, j+1\}$ is equal to $\tau e$, so by (\ref{(1.12)}),
 $\L g$ is discontinuous at
$\tau e$.

Now consider the remaining case where $j \in \{0,q\}$. We treat the case $j = q$; the other case is similar. Since $g$ is discontinuous
at
$e = x_j$, then
$n_j
\not= 0$, so
$n_j\tau([x_{j-1},x_j])$ is discontinuous at $\tau x_j = \tau e$. As in the previous paragraph,  none of the
endpoints of the intervals
$\tau([x_{i-1},x_{i}])$ for $i\not= j$ is equal to $\tau e$, so by (\ref{(1.12)}),
 $\L g$ is discontinuous at
$\tau e$.}

\lem{1.46}{Let $\tau:I\to I$ be \pwm,  with associated partition $C$, and let $C_0 = C\setminus \{0,1\}$. Let $g
\in C(X,\Z)$, and $e \in \D(g)\setminus C_0$.  If for all $n \ge 0$,
\begin{enumerate}
\item  $\tau^n e \notin \htau^k C_0$ for $k \le n$, and\label{1.46i}
\item $\tau^n e \notin \htau^n (\D(g)\setminus \{e\})$,\label{1.46ii}
\end{enumerate}
\noindent then for all $n \ge 0$, $\L^ng$ is discontinuous at $\tau^n e$.}

\prooff{By (\ref{(1.13)}),
\begin{equation}
\D(\L^n g) \subset \htau^n\bigl(\D(g)\setminus \{e\}\bigr)\cup \{\tau^ne\} \cup \bigcup_{k=1}^n
\htau^k C_0.
\label{(1.14)}
\end{equation}
We will prove
\begin{equation}
\tau^ne \in \D(\L^n g)\label{(1.15)}
\end{equation}
for all $n$, by induction on $n$. Suppose (\ref{(1.15)}) holds for a particular $n\ge 0$. We will apply Lemma \ref{1.45} with
$\tau^ne$ in place of $e$. By (\ref{1.46i}), $\tau(\tau^ne) \notin \htau C_0$. Suppose $x \in
\tau^{-1}\tau(\tau^ne) \cap \D(\L^n g)$, so that $\tau x = \tau^{n+1}e$. We will show that $x$ is not in the
first or third set on the right side of (\ref{(1.14)}), so must be in the second. By (\ref{1.46ii}), 
$\tau x =
\tau^{n+1}e \notin \htau^{n+1}(\D(g)\setminus \{e\})$, so $x \notin \htau^n(\D(g)\setminus \{e\})$. Similarly,
if $x \in \htau^k C_0$ with
$k
\le n$, then $\tau^{n+1}e = \tau x \in \htau^{k+1}C_0$, which contradicts (\ref{1.46i}).  Thus
(\ref{(1.14)}) implies $x =
\tau^ne$, so we've shown
\begin{equation}
\tau^{-1}\tau(\tau^ne) \cap \D(\L^ng) = \{\tau^ne\}.
\end{equation}
By Lemma \ref{1.45}, $\tau^{n+1}e \in \D(\L^{n+1}g)$, so by induction, $\tau^ne \in \D(\L^ng)$ for all $n \ge 0$.
}

Recall that a $\Z[t,t^{-1}]$ module $G$ is {\it cyclic\/} if there is an element $e \in G$ such that
$G = \Z[t,t^{-1}]e$.

\prop{1.47}{Let $\tau:I\to I$ be \pwm,  with associated partition $C$. Assume that $DG(\tau)$ is cyclic,
and that there exists $a \in \{0,1\}$ with an infinite orbit, such that 
\begin{equation}
\htau (C\setminus \{a\}) \subset C.\label{(1.16)}
\end{equation}
 Then $DG(\tau)\cong \Z[t,t^{-1}]$ as abelian groups, with the action of $\L_*$ given by
multiplication by $t$.}

\prooff{Let $b \in C$ be adjacent to $a$ (i.e., there are no points in $C$ between $a$ and $b$). Since the
orbit of $a$ is infinite, $\tau^n a \not= \tau^k a$ for $k \not= n$, and $\tau^n a \notin C$ for $n \ge 1$.
By (\ref{(1.16)}) and Lemma \ref{1.46} (with
$e = a$ and
$g = I(a,b)$), for each
$n
\ge 0$,
$\L^n I(a,b)$ is discontinuous at
$\tau^n a$.  On the other hand, by (\ref{(1.13)}), for $k < n$,
\begin{equation}
\D(\L^k I(a,b)) \subset \{\tau^k a\}\,\cup\, \htau^k b\, \cup\, \bigcup_{j=1}^k \htau^j C_0,
\end{equation}
where $C_0 = C\setminus \{0,1\}$. It follows that $\L^k I(a,b)$ is continuous at $\tau^n a$ for $k < n$.
Therefore, for any polynomial $p\in \Z[t]$ of degree $n$, $p(\L)I(a,b)$ is
discontinuous at $\tau^n a$, and in particular is not identically zero. 

By hypothesis, $DG(\tau)$ is cyclic, so there exists $g\in C(X,\Z)$ such that $[g]$ is a
generator for
$DG(\tau)$. Suppose that
$q(\L_*)[g] = 0$ for some $q\in \Z[t,t^{-1}]$. Choose $q'\in \Z[t,t^{-1}]$ so that $q'(\L_*)[g] =
I(a,b)$. Then $q(\L_*)q'(\L_*)I(a,b) = 0$. It follows that there is a non-zero
polynomial
$q''\in Z[t]$ such that $q''(\L)I(a,b) = 0$, contrary to the result established in the first
paragraph.  We conclude that  $p \to p(\L_*)[g]$ is an
isomorphism from $\Z[t,t^{-1}]$ onto $DG(\tau)$.}

We will apply Proposition \ref{1.46} to find dimension groups of unimodal maps 
later in this paper, cf. Theorem \ref{1.88}, and we will show for multimodal maps that a similar requirement of 
infinite disjoint orbits for critical points leads to the dimension module being free, cf. Proposition \ref{1.53}.

\section{Markov maps}

\defi{1.80}{A \pwm\ map $\tau$ is {\it Markov\/} if there is
a partition $0 = b_0 < b_1 < \ldots < b_n= 1$, with each $b_i$ being in $I_1$,
such that for each $i$, $\tau$ is monotonic on $(b_i,
b_{i+1})$,  $\overline{\tau(b_i,b_{i+1})}$ is 
 a union of intervals of the form
$[b_j,b_{j+1}]$, and  for some $k \ge 0$, $\htau^k(C) \subset \{b_0, b_1, \ldots, b_n\}$.  We
call such a partition a {\it Markov partition for $\tau$\/}, and will refer to $(b_0, b_1)$,
$(b_1,b_2)$, $\ldots$, $(b_{n-1}, b_n)$ as the {\it partition intervals for $\tau$}.}

For a Markov map, we have both the original partition associated with $\tau$, and the Markov
partition. Since we don't require that
$\tau$ be continuous on each interval
$(b_{i-1},b_i)$, a Markov partition for $\tau$ may not even qualify as a ``partition associated with
$\tau$" in the sense defined at the start of this paper. 

If $C$ is the set of endpoints in
the  partition associated with 
$\tau$, then $\tau$ will be Markov iff the forward orbit of $C$
under
$\htau$ is finite.  In that case,  the points of the
forward orbit of $C$ will be the endpoints of a Markov partition.  When
$\tau$ is continuous, $\tau$ will be Markov iff each point in $C$
is eventually periodic. 

If $\tau:I\to I$ is \pwm, and $\sigma:X\to X$ is the associated
local homeo\-morphism, a {\it Markov partition for $\sigma$\/} is a
partition of $X$ into clopen order intervals $E_1, \ldots E_k$
such that $\sigma$ is monotonic on each $E_i$, and maps $E_i$
onto a union of some of $E_1, \ldots, E_k$. If  $0 = b_0 < b_1 <
\ldots < b_n= 1$ is a Markov partition for $\tau$, then
defining $E_i = I(b_{i-1},b_i)$ gives a 
Markov partition for $\sigma$.  
The associated \zeroOne {\it incidence
matrix\/} $A$ is given by $A_{ij} = 1$ iff $\sigma(E_i) \supset
E_j$, or equivalently, if $\tau(b_{i-1},b_i) \supset (b_{j-1},b_j)$..

Note that if $E_1, \ldots, E_n$ is a Markov partition for
$\sigma$, then 
$\{\sigma^k(X)\}$ will be a nested decreasing sequence of sets,
with each being a union of some of $E_1, \ldots, E_n$.  Thus
$\sigma$ will be eventually surjective (Definition \ref{1.23}), and similarly so will $
\tau$ (Lemma \ref{1.23.1}).

\defi{1.81}{If $A$ is an $n \times n$ \zeroOne  matrix, $G_A$ denotes the stationary inductive limit
 $\Z^n \mapright{A}\Z^n$ in the category of ordered abelian groups. Here $A$ acts by right multiplication,
and 
$G_A$ will be a dimension group.  The action of $A$ induces an automorphism of $G_A$ denoted $A_*$, so we
view
$G_A$ as a
$\Z[t,t^{-1}]$ module, and refer to $(G_A,G_A^+, A_*)$ as the dimension triple associated with $A$.}

\notate{Let $G$ be an ordered abelian group, $T:G\to G$ a positive homomorphism, and $G_1, G_2, \ldots $  a sequence of
subgroups of
$G$ such that $T(G_i) \subset G_{i+1}$. If $g \in G_n$, we will write $(g,n)$ for the sequence $(0, 0, \ldots, 0, g, Tg, T^2g,
\ldots)$, where $n-1$ zeros precede $g$. We write $[(g,n)]$ (or simply $[g,n]$) for the equivalence class of this sequence in
the inductive limit. (See the notation introduced after Definition \ref{1.16}).  Every element of the inductive limit is of the
form
$[g,n]$, and
$[g_1,n_1] = [g_2,n_2]$ with
$n_2
\ge n_1$ iff $T^{n_2-n_1+k}g_1 = T^kg_2$ for some $k \ge 0$, or equivalently, iff $T^{n_2+k}g_1 = T^{n_1+k}g_2$. In particular,
we represent elements of
 $G_A$ as pairs $[v,k]$, where $v \in \Z^n$ and $k \in \N$, and the automorphism $A_*$ of $G_A$ is given by $A_*([v,k]) =
[vA,k]$. Similarly, there is a canonical automorphism associated with any stationary inductive limit.  For an exposition of
inductive limits for dimension groups, see \cite{Ror}.}
 
\lem{1.81.1}{Let $C_1, C_2, \ldots$ be subgroups of $C(X,\Z)$ such
that
$\L C_n \subset C_{n+1}$ for all $n$, and let $G$ be the
inductive limit of the sequence $C_n \mapright{\L}C_{n+1}$.  If
\begin{enumerate}
\item  for each
$f
\in C(X,\Z)$ there exists $n \ge 0$ such that $\L^n f \in
C_n$,\label{1.81.1i} 
\item for every $k \ge 0$ and $f \in C_k$, there exists $g \in C(X,\Z)$ such that $\L_*^k[g]=
[f]$,\label{1.81.1ii}
\end{enumerate}
then the map $\Phi:G\to
DG(\tau)$ defined by
$\Phi([f,n]) = \L_*^{-n}[f]$ is an order isomorphism from the dimension group
$G$ onto the dimension group $DG(\tau)$.  If $C_1 = C_2= C_3 \cdots$, then this isomorphism carries
the canonical automorphism of this stationary inductive limit to $\L_*$.} 

\prooff{Recall that $\L_*$ is injective, but need not be surjective. However, by (\ref{1.81.1ii}),
$\L_*^{-n}$ makes sense on
$C_n$. Suppose that
$f_1
\in C_{n_1}$ and
$f_2
\in C_{n_2}$, and that
$[f_1,n_1] = [f_2,n_2]$. Then $\L^{n_2+k}f_1 = \L^{n_1+k}f_2 \in C_{n_1+n_2+k}$ for some $k \ge 0$,
so taking equivalence classes and applying $\L_*^{-n_1-n_2-k}$
gives $\L_*^{-n_1}[f_1] = \L_*^{-n_2}[f_2]$. Thus $\Phi$ is
well defined. 

If $f_1 \in C_{n_1}$, $f_2 \in C_{n_2}$ and
$\Phi([f_1,n_1]) = \Phi([f_2, n_2])$, then reversing the
previous argument gives $\L_*^{n_2}[f_1]  =
\L_*^{n_1}[f_2]$, so for some $n \ge 0$, $\L^{n+n_2}f_1 =
\L^{n+n_1}f_2 \in C_{n+n_1+n_2}$. Thus 
$[f_1,n_1] = [f_2,n_2]$, proving that $\Phi$ is 1-1.

For any $f \in C(X,\Z)$, choose $n$ so that $\L^n f\in C_n$.
Then $\Phi[\L^n f, n] = \L_*^{-n} [\L^n f] = [f]$, so $\Phi$
is surjective.  The remaining assertions are readily verified.}

\prop{1.82}{Let  $\tau:I\to I$ be \pwm\ and Markov, with associated local homeomorphism $\sigma:X\to X$.
Let
$E_1, E_2,
\ldots, E_q$ be the associated Markov partition for $\sigma$, with incidence matrix
$A$, and define $\psi:\Z^q \to C(X,\Z)$ by $\psi(z_1, z_2, \ldots, z_{q}) = \sum_i z_i E_i$.  Then the map 
$\Phi:G_A\to DG(\tau)$ defined by
$\Phi([v,n]) = \L_*^{-n}[\psi(v)]$ is an isomorphism from the dimension triple
$(G_A,G_A^+,A_*)$ onto the dimension triple $(DG(\tau),DG(\tau)^+,\L_*)$.}

\prooff{Let $M$ be the range of $\psi$.  Then $\psi$ is an order
isomorphism from $\Z^q$ onto $M$, and $M$ is invariant under $\L$ (by the Markov
property of the partition).  Observe that 
\begin{equation}
\L(\psi(v))
= \psi(vA) \hbox{ for all $v \in \Z^q$},\label{(1.33)}
\end{equation}
so $\psi$ is an isomorphism from $(G_A,G_A^+, A_*)$ onto the inductive limit of the stationary sequence
$\L:M\to M$.

 We next show that for any $f \in C(X,\Z)$ there exists
$n \ge 0$ such that $\L^n f \in M$. It suffices to prove this
for $f = I(a,b)$ with $a, b \in I_1$. Note that $\L I(a,b)$ is
a sum of intervals $I(c,d)$ with $c, d$ either being images of
$a, b$ or being points in $\htau C$.  By definition of $I_1$, for each
point
$x \in I_1$ there exists $n$ such that $\tau^n(x) $ is in the
forward orbit of $C$. It follows that for $k$ large enough, $\L^k
f$ is  a sum of intervals of the form $I(x,y)$ with endpoints $x, y$ in the forward orbit of
$C$, and thus, by the definition of a Markov partition (Definition \ref{1.80}), 
for $k$ large enough
these endpoints will be in
$\{b_0, \ldots, b_q\}$. Thus  $\L^n f \in M$ for $n$ large enough. Since $\tau$ is Markov, $\tau$ and $\sigma$ are eventually
surjective, so $\L_*$ is surjective (Proposition \ref{1.26}). Now the conclusion follows from Lemma \ref{1.81.1} with each of
$C_1$,
$C_2$,
$\ldots$ equal to
$M$.}

Every \zeroOne square matrix (without zero rows) is the incidence matrix of a Markov map. To illustrate,  the matrix 
\begin{equation}
A =
\begin{pmatrix}
0&1&1\cr 1&0&1\cr 
1&1&0
\end{pmatrix}\label{(1.34)}
\end{equation}
is the incidence matrix for the Markov map pictured in Figure \ref{fig1} on page~\pageref{fig1},
with respect to the Markov partition $B = \{0, 1/3, 2/3, 1\}$. In this example, the maximal partition
associated with $\tau$ is $C = \{0, 1/3, 1/2, 2/3, 1\}$.

Let $A$ be an $n \times n$ \zeroOne matrix (with no zero rows), and let $(X_A, \sigma_A)$ be 
the associated one-sided shift of finite type. (In other words, $X_A$ consists of sequences
$x_0x_1x_2\ldots$ with entries in $\{1,2,\ldots, n\}$ such that $A_{x_ix_{i+1}} =
1$ for all
$i$, and
$\sigma$ is the left shift.)  Then,
for an appropriate metric,
$X_A$ is a zero dimensional metric space, and $\sigma_A$ is a local homeo\-morphism. The dimension group
$G_{\sigma_A}$ (cf. Definition \ref{1.12}) is isomorphic to
$G_A$.  (For
$A$ irreducible, this is \cite[Thm. 4.5]{BFF}; the proof in \cite{BFF} works in general. Alternatively,  the proof of
Proposition \ref{1.82} can be adapted to prove this.) This 
suggests that $(X,\sigma)$ might be conjugate to the one-sided shift
$(X_A,\sigma_A)$, and we now show that this happens precisely if $\tau$ has no
homtervals (Definition \ref{1.36}). 

Given a Markov partition $E_1, \ldots, E_n$ for $ \sigma:X\to X$, with incidence matrix $A$,
the {\it itinerary map $S:X\to X_A$} is given by $S(x) = s_0s_1s_2\ldots$, where
$\sigma^k(x) \in E_{s_k}$. We say \emph{itineraries separate points of $X$} if the itinerary map
is 1-1.  The following proposition shows that this property is independent of the choice of Markov
partition.

\prop{1.83}{Assume $\tau:I\to I$ is \pwm\ and Markov, with incidence
matrix
$A$. These are equivalent.
\begin{enumerate} 
\item $\tau$ has no homtervals.\label{1.83i}
\item Itineraries separate points of $X$.\label{1.83ii}
\item $(X,\sigma)$ is conjugate to the one-sided shift of finite type
$(X_A,\sigma_A)$.\label{1.83iii}
\end{enumerate}
In particular, these equivalent conditions hold if $\tau$ is transitive.}

\prooff{ Let $B = \{b_0, b_1 ,\ldots, b_q\}$ be the given Markov partition for $\tau$, and let
$C$ be the partition associated with $\tau$.

(\ref{1.83i}) $\Rightarrow$ (\ref{1.83ii}) Assume that there are no homtervals. We  will show that
the itinerary map
$S:X\to X_A$ is 1-1.   Suppose (to reach a contradiction) that there exists points  $x< y$ in
$X$ with the same itinerary. 

 We first consider the case where $x = a^-$ and $ y = a^+$ for $a \in I_1$.  By definition of $I_1$, the
orbit of
$a$ eventually lands on some point in the forward orbit of
$C$, and by the definition of a Markov partition, the forward orbit of $C$ is  eventually contained in $B$. 
Thus there exists $n$ such that $\htau^n a \subset B$, so $\pi(\sigma^n(a^\pm)) \subset B$.  Suppose that
$\sigma^k(a^-)$ and
$\sigma^k(a^+)$ belong to the same Markov partition interval for $0\le k\le n-1$. Then
$\sigma$ is monotonic on the Markov partition interval containing $\sigma^{n-1}(a^-)$ and $\sigma^{n-1}(a^+)$, so
$\sigma^n(a^-) \not=
\sigma^n(a^+)$. If $\sigma^n(a^-) = b^\pm$ for $b \in B$, then $\sigma^n(a^+) = b^\mp$, so
the itineraries of $a^\pm$ are different, contrary to our assumption.
 
Now assume 
$(x,y)_X$ is not empty, and that $x, y$ have the same itinerary. By density
of $X_0$ in $X$, we may assume $x$ and $y$ are in $X_0$.  Then $x'=\pi(x)$ and $y'=\pi(y)$ are distinct
points in
$I_0$ with the same itinerary with respect to the given Markov partition of $I$. Let $J= (x',y')$, and for $n \ge 0$ define $J_n$ to be
the open interval whose endpoints are
$\tau^n(x')$ and $\tau^n(y')$.  Since $x$ and $y$ have the same itinerary, $\tau$ is  monotonic on $J_n$ for all $n$.
 Since
there are no homtervals, there is some
$n$ such that $J_n$ contains a point $c \in C$.
Then  $\htau^k(c) \subset J_{n+k}$ for
all $k \ge  0$.  However,  the orbit of $c$ eventually lands in $B$, so $J_q$ contains
a point of $B$ for some $q$. But then $\tau^q(x')$ and $\tau^q(y')$ would have different
itineraries, a contradiction. Thus $S$ is 1-1.

(\ref{1.83ii}) $\Rightarrow$ (\ref{1.83iii}) Assume itineraries separate points, so that $S$ is 1-1.
To show that $S$ is surjective, let $s = s_0s_1\cdots
\in
X_A$. For each $n \ge 0$, let
$I_{s_0s_1\ldots s_n}$ be the set of points in $X$ with
the initial itinerary $s_0s_1\ldots s_n$.  Then the points with itinerary $s$ are those that
are in $\cap_{n=0}^\infty I_{s_0s_1\ldots
s_n}$. The sets $I_{s_0s_1\ldots
s_n}$ for $n = 0, 1, \ldots$ are nested and compact, so surjectivity of the map $S$ will follow if we show
each such set is nonempty.

Observe that $I_{s_0s_1\ldots
s_n}=
I_{s_0}
\cap \sigma^{-1}(I_{s_1\ldots s_n})$.   We use induction on the length of the partial
itinerary. Suppose that these partial itinerary sets are nonempty  for partial itineraries of length $m$. If
$s_0s_1$ is an allowed transition with respect to the incidence matrix $A$, by definition of $A$,  $\sigma(I_{s_0}) \supset I_{s_1}$, so
$\sigma(I_{s_0}) \supset I_{s_1\ldots s_m}$. Thus
$\sigma(I_{s_0})\cap\sigma^{-1}(I_{s_1\ldots s_m})$ is not empty, which completes the proof that $S$
is surjective, so that $(X,\sigma)$ is conjugate to $(X_A, \sigma_A)$.

(\ref{1.83iii}) $\Rightarrow$ (\ref{1.83i}) If $(X,\sigma)$ is conjugate to a shift of finite type,
then $X$ is totally disconnected. 
 By Proposition
\ref{1.37},   $\tau$ has no homtervals.

 Finally, if $\tau$ is transitive, then $\tau$ has no homtervals (Proposition
\ref{1.37}).}

\exe{1.83.1}{Let $\tau_0(x) = k x(1-x)$, with $c$ the unique critical point, and choose $k \approx 3.68$  so that $p =
\tau^3(c)$ is fixed. Let $\tau$ denote $\tau_0|_J$, where  $J = [\tau^2(c),
\tau(c)]$, rescaled so that the domain of the restricted map is $[0,1]$. (See the restricted logistic map in Figure \ref{fig1} on
page~\pageref{fig1}.) The map
$\tau$ is Markov, with Markov partition $\{0,c,p,1\}$. The incidence matrix is 
$$A = \begin{pmatrix}0&0&1\\ 0 &0&1\\ 1&1&0\end{pmatrix}.$$
Since $\tau$ is a polynomial, it will have homtervals iff it has an attracting periodic
orbit (Proposition \ref{1.36.1}). For a quadratic, the orbit of the critical point will be attracted to any attracting
periodic orbit (\cite[p. 158]{Devaney}).  For our particular polynomial $\tau_0$, since the orbit of the critical point $c$ lands
on a repelling fixed point
$p$, we conclude that there are no homtervals. If $\sigma:X\to X$ is the associated local homeo\-morphism, by Proposition \ref{1.83},
$(X,\sigma)\cong (X_A, \sigma_A)$.}

If $\tau:I\to I$ is \pwm\ and Markov, the {\it directed graph associated with $\tau$\/} has vertices consisting of the
Markov partition intervals, with an edge from $E$ to $F$ if $\tau(E) \supset F$. This is the same as the directed graph
associated with the incidence matrix for $\tau$. By ``directed graph" we will always mean a graph in which there are no
multiple edges between pairs of vertices.  A \emph{loop} in a directed graph is a path that begins and ends at the same vertex. If
$\tau:I\to I$ is piecewise linear and Markov, we will see we can determine whether  there are homtervals by looking at the directed
graph associated with
$\tau$.

\defi{1.84}{Let  $v_1, \ldots, v_n$ be
vertices of a directed graph. If $v_1v_2\ldots v_n v_1$ is a loop, then an {\it
exit\/} for this loop is an edge from some
$v_i$ to a vertex other than
$v_{i+1\bmod n}$.  A directed graph satisfies \emph{Condition L} if every loop has an exit.}

Condition L appears  in the
study of C*-algebras associated with (possibly infinite) directed graphs, cf. \cite{KumPaskRae}.  An equivalent  condition
 for finite graphs (``Condition I") was assumed in
\cite{CK}. If a (finite) directed graph is irreducible,  Condition L is equivalent to the
associated incidence matrix not being a permutation matrix. 

In Proposition \ref{1.85}, by a piecewise linear Markov map we mean a \pwm\ map, which is Markov, and which on
each interval of the Markov partition is piecewise linear. We allow the slopes of the pieces to vary from interval to
interval, but not within an interval of the partition.  We don't require the map to be continuous on each partition
interval, since we want to  allow examples like the Markov map in Figure \ref{fig1} on page~\pageref{fig1}.

\prop{1.85}{Let $\tau:I\to I$ be  piecewise linear and Markov, with incidence matrix $A$. These are equivalent. 
\begin{enumerate} 
\item $\tau$ has no homtervals.\label{1.85i}
\item Itineraries  separate points of $X$.\label{1.85ii}
\item $(X,\sigma)$ is conjugate to the one-sided shift of finite type
$(X_A,\sigma_A)$.\label{1.85iii}
\item In the directed graph associated with $\tau$, every loop has an exit.\label{1.85iv}
\end{enumerate}
}

\prooff{By Proposition \ref{1.83}, it suffices to prove that (\ref{1.85iv}) and (\ref{1.85i}) are
equivalent. 
 Suppose $v_1v_2\ldots v_nv_1$ is a loop without an exit in the directed graph associated with $\tau$. By definition, each $v_i$ is
an open interval in the associated Markov partition for $\tau$. Since $\tau$ is monotonic and piecewise linear on each $v_i$, and
$\overline{\tau(v_i)} = \overline{v_{i+1\bmod 1}}$ for each $i$, then $\tau$ must be a homomorphism of $v_i$ onto $v_{i+1\bmod n}$ for
each $i$. Thus $v_i $ is a homterval.

Conversely, assume that every loop has an exit.  By conjugating $\tau$ by a suitable continuous
piecewise linear map if necessary, we may assume that each interval in the Markov partition has the same length.
If $E$ is an interval in the Markov partition, since  the slope of $\tau$ within $E$ is constant, it must be
an integer, equal to the number of partition intervals covered by $\tau(E)$.  Let
$V$ be a homterval contained in a partition interval $E_0$. We claim that there exists $n \ge 1$ such that the length
$\tau^n(V)$ is at least double that of $V$.  If the slope of $\tau|_{E_0}$ is $\pm 1$, since all partition
intervals have the same length, then
$E_1 = \tau(E_0)$ will also be a partition member.  If $E_0, \ldots, E_j$ have been defined, and the slope of
$\tau$ on $E_j$ is $\pm 1$, define $E_{j+1} = \tau(E_j)$. Note that if $E_i = E_j$ for $i < j$, then $E_i, \ldots,
E_j$ would be a loop without an exit, which is not allowed, so the sets $E_1, E_2, \ldots, E_j$ are
distinct. Choose $n$ maximal such that $E_0$, $E_2 = \tau(E_1)$, $\ldots$, $E_n = \tau(E_{n-1})$ are all
partition members. Then $\tau|_{E_n}$ has slope $\ge 2$, so the (Lebesgue) measure of $\tau^n(V)$ is at least double the
length of
$V$. Furthermore, since $V$ is a  homterval, $\tau$ is 1-1 on $\tau^k(V)$ for all $k$, so the measure of $\tau^{k}(V)$ never
decreases as $k$ increases. 

Now let $V$ be an arbitrary homterval. Then we can write $V$ as a disjoint union of homtervals $V_1, \ldots, V_k$,
each contained in a partition interval, together with a finite set of points. Then there exists $n$ such that 
for all $i$, $\tau^n (V_i)$ has measure at least double the length of $V_i$. Then the measure of $\tau^n(V)$ is at least
double the length of
$V$.  This process can be repeated indefinitely, which is impossible since $\tau^n(V) \subset [0,1]$ for all $n$. We
conclude that $\tau$ has no homtervals.}

This provides an easy way to tell if Markov piecewise linear
maps are transitive or topologically exact.  Recall that a \zeroOne matrix is {\it irreducible} if in the
associated directed graph, there is a path from any vertex to any other vertex, and is  {\it
irreducible and aperiodic\/} if it is irreducible and the greatest common divisor of the lengths of paths between any two
vertices is 1.  This is equivalent to $A$ being {\it primitive\/}, i.e., to some power of $A$ having strictly positive entries. It is
well known that the one-sided shift $\sigma_A$ is transitive iff $A$ is irreducible, (e.g., \cite[Thm. 1.4.1]{Kitchens}), and is
topologically mixing iff $A$ is primitive. It is readily verified that $\sigma_A$ will be topologically exact iff it is
topologically mixing.

\cor{1.86}{Assume $\tau:I\to I$ is piecewise linear and Markov,  with
incidence matrix
$A$.  Then  $\tau$ is topologically exact iff $A$ is primitive, and is transitive iff $A$ is irreducible and 
is not a permutation matrix. If $A$ is irreducible, and is not a permutation matrix, and  $\sigma:X\to X$ is the  local homeo\-morphism
associated with $\tau$, then $(X,\sigma)$ is conjugate to the one-sided shift $(X_A,\sigma_A)$.}

\prooff{Assume $\tau$ is transitive. By Proposition \ref{1.83}, $(X,\sigma)$ is conjugate to the one-sided shift
$(X_A,\sigma_A)$. Transitivity of $\tau$ implies transitivity of $\sigma$ (Proposition \ref{1.8}). Thus $\sigma_A$
is transitive, so $A$ is irreducible. If $A$ were a permutation matrix, then $A^p$ would be the identity for some $p$, so 
$\tau^p$ restricted to each  interval of the Markov partition would be a monotonic bijection of that interval onto itself. This
contradicts transitivity of $\tau$, so
$A$ cannot be a permutation matrix.

 Conversely, if $A$ is irreducible and not a permutation matrix, then $\sigma_A$ is  transitive, and every loop has an exit, which 
implies that
$\tau$ has no homtervals (Proposition \ref{1.85}).   Having no homtervals implies that
 $(X,\sigma)$ is conjugate to $(X_A, \sigma_A)$ (Proposition \ref{1.83}). Thus $\sigma$ and $\tau$ are transitive.

If $\tau$ is topologically exact, then it is also transitive, so again $(X,\sigma)$ is
conjugate to $(X_A, \sigma_A)$. Then $\sigma$ and $\sigma_A$ are topologically exact, so  
$A$ is primitive.   Conversely, if $A$ is primitive, one easily verifies that in the associated graph every loop has an exit, so
$(X,\sigma)$ will be conjugate to $(X_A, \sigma_A)$. Primitivity of $A$ implies that $\sigma_A$ is topologically
exact, and thus that $\sigma$ and $\tau$ are topologically exact.}

\exe{1.87}{Let $\tau$ be the Markov map in Figure \ref{fig1} on page~\pageref{fig1}, and let $\sigma:X\to X$ be the associated local
homeo\-morphism. Then $\tau$ is Markov, so the dimension triples
$(DG(\tau),DG(\tau)^+, \L_*)$ and $(G_A, G_A^+, A_*)$ are isomorphic (Proposition \ref{1.82}). The
associated incidence matrix
$A$ is primitive, so $(X,\sigma)$ is conjugate to $(X_A, \sigma_A)$ and $\tau$ is topologically exact (Corollary \ref{1.86}).}

\section{Unimodal maps}

If $\tau:I\to I$ is continuous, with just two intervals of monotonicity,
\label{unimodal} without loss of generality, we may assume $\tau$ increases and then decreases.  (If not,
conjugate by the map $\phi(x) = 1-x$.)  We will say a continuous map $\tau:I\to I$ is
{\it unimodal\/} if there exists $c$ in
$(0,1)$ such that $\tau$ is increasing on $[0,c]$ and decreasing
 on $[c,1]$. We are mainly interested in the case when $\tau$ is eventually surjective.
In that case, for the purpose of computing the dimension group, we may as well assume that $\tau$ is
surjective, cf. Proposition \ref{1.28}.

\theo{1.88}{Let $\tau:I\to I$ be unimodal and surjective. Then the dimension
module
$DG(\tau)$ is cyclic, with generator $I(0,1)$. If the critical point
$c$ is eventually periodic, then $\tau$ is Markov.  In that case, if the incidence matrix
is $A$, then $(DG(\tau), DG(\tau)^+, \L_*) \cong (G_A, G_A^+, A_*)$.  If $c$ is not
eventually periodic, then
$DG(\tau)
\cong \Z[t,t^{-1}]$ (as abelian groups), with the action of $\L_*$ given by multiplication by $t$.}

\prooff{Suppose $\tau(1) = 0$. By Corollary \ref{1.42}, the module $DG(\tau)$ is generated by
$I(0,c)$ and
$I(c,1)$, where $c$ is the unique critical point. Since
$\L_* I(c,1) = I(0,1)$, then $I(0,c) = I(0,1) - (\L_*)^{-1}I(0,1)$, so $DG(\tau)$  is
generated by $I(0,1)$, and thus is cyclic. The same conclusion applies if $\tau(0) = 0$,
by a similar argument.

If $c$ is eventually periodic, then $\tau$ is Markov with respect to the partition given
by the forward orbit of $c$, together with the endpoint 0. Then $(DG(\tau), DG(\tau)^+, \L_*) \cong (G_A, G_A^+, A_*)$
 follows from Proposition \ref{1.82}. If
$c$ is not eventually periodic, by surjectivity the unique critical point maps to 1,  one
endpoint of $I$ must map to zero, and the other endpoint will have an infinite orbit. Thus $DG(\tau) \cong \Z[t,t^{-1}]$ by
Proposition \ref{1.47}.}

Given $s> 1$, we define the {\it restricted tent map\/} $T_s$ by 
\begin{equation}
T_s(x) = \begin{cases}
1+s(x-c) & \text{if $x \le c$}\cr
1-s(x-c) &\text{if $x > c$},
\end{cases}
\end{equation}
where $c = 1-1/s$. 
(This is the usual symmetric tent map on $[0,1]$ with slopes $\pm s$, restricted to the interval $[\tau^2(1/2), \tau
(1/2)]$, which is the interval of most interest for the dynamics.    Then the map has been rescaled so that its domain
is  [0,1]. See Figure \ref{fig1} on page~\pageref{fig1}.) Note that $\tau(c) = 1$ and $\tau(1) = 0$.
\smallskip

For use in the next example, we recall one way to compute the dimension group $G_A$ (Definition \ref{1.81}), due to Handelman;
see
\cite{Lind-M} for an exposition. If
$A$ is
$n
\times n$, view
$A$ as acting on $\Q^n$, and let $V_A\subset \Q^n$ be the eventual range of $A$, i.e.
$V_A = \cap_{k=1}^\infty
\Q^n A^k$.  Let
$G$ be the set of vectors $v \in V_A$ such that $vA^k \in \Z^n$ for some $k \ge 0$, and $G^+$ those
vectors such that $vA^k \in (\Z^n)^+$ for some $k \ge 0$. Then $(G_A, G_A^+,A_*) \cong (G, G^+, A)$.

\exe{1.92}{Let $\tau = T_2$ be the full tent map.  Then $\tau$ is Markov,
with incidence matrix
$A = \begin{pmatrix} 1&1\cr1&1\end{pmatrix}$, so $DG(\tau) \cong G_A$. Since the incidence matrix is primitive, $\tau$ is
topologically exact (Corollary \ref{1.86}), as also is easy to verify directly. 
We compute
$G_A$. The eventual range of $A$ in $\Q^2$ consists of multiples of
$(1,1)$, and multiplication by $A$ on the eventual range multiplies this vector by 2. Thus $G_A$ is isomorphic to the space of
vectors $q(1,1)$ for $q$ in $\Q$ such that $q(1,1)A^k = 2^k q(1,1) \in \Z^2$ for some $k \ge 0$. 
It follows that $DG(\tau)$ is isomorphic to the group of dyadic
rationals, with the  automorphism $\L_*$ given by multiplication by 2.}

\exe{1.92.1}{Let $\tau$ be the restricted tent map $T_{\sqrt{2}}$,  let $c$ be the critical point, and let $p$ be the fixed
point.  Then
$\tau$ is  Markov, with Markov partition
$\{0,c,p,1\}$, and $\tau^2$ is conjugate to the tent map on each of the invariant intervals $[0,p]$ and
$[p,1]$. Thus $\tau$ is transitive (as can also be seen by applying Corollary \ref{1.86}).  By Example
\ref{1.92},   $DG(\tau)$ will be the direct sum of two copies of the dyadic
rationals, with the generators of the summands being $I(0,p)$ and $I(p,1)$. The associated automorphism $\L_*$ will exchange the two
summands, taking
$(a,b) \in \Z[1/2] \oplus \Z[1/2]$ to
$(b,2a)$.}

\exe{1.94}{Let $s = 3/2$. Then the critical point of the restricted tent map $T_s$  is not
eventually periodic. (Indeed, if $c$ is the critical point and  $n \ge 3$, then $\tau^n (c)  = a_n/2^{n-2}$, where $a_n$
is odd.) By Theorem \ref{1.88},
$DG(T_s)
\cong
\Z[t,t^{-1}]$ as abelian groups, with the action of $\L_*$ given by multiplication by $t$.} 

\exe{1.95}{Let $\tau(x) = k x(1-x)$, with $k$ as in Example \ref{1.83.1}, restricted to the interval described there. 
Then 
$\tau$ is Markov, with incidence matrix is 
$$A = \begin{pmatrix}0&0&1\\ 0 &0&1\\ 1&1&0\end{pmatrix}.$$
As shown in Example \ref{1.83.1}, the associated local homeo\-morphism $(X,\sigma)$ is conjugate to
$(X_A,\sigma_A)$.  The map $T_{\sqrt{2}}$ is Markov, with the same
incidence matrix. Since $A$ is irreducible and not a permutation matrix, 
$T_{\sqrt{2}}$ also is  conjugate to $(X_A,\sigma_A)$
(Corollary \ref{1.86}). Thus $\tau$ is conjugate to $T_{\sqrt{2}}$, so $DG(\tau)$ is as described in Example \ref{1.92.1}. 
Alternatively, $DG(\tau)$ can be computed by applying Proposition \ref{1.82}.}

\section{Multimodal maps}

By a multimodal map we mean a continuous \pwm\ map. In this section we will show that the dimension
modules for multimodal maps are free modules, under the assumption that the critical points have orbits that
are infinite and disjoint.
 
\defi{1.48}{Let $\tau:I\to I$. If a set $B \subset I$ has the property that for 
$b_1, b_2\in B$, and $n, m \ge 0$,
\begin{equation}
\tau^m b_1 = \tau^n b_2\implies b_1 = b_2 \hbox{ and } m = n,
\end{equation}
we say $B$ satisfies the infinite  disjoint orbit condition (IDOC).}

Note that if $B$ satisfies the IDOC, then $\tau$ is 1-1 on $\orbit(B)$. The IDOC
condition  has been used in studies of interval exchange maps,  cf., e.g., \cite{KeaneDisc}.

Recall from Section \ref{Cyclic} that we view functions in $C(X,\Z)$ as functions in $L^1(\R)$, continuous except on a finite
set, and that
$\D(g)
\subset I_1$ denotes the set of points where $g \in L^1(\R)$ has an essential (jump) discontinuity.

\lem{1.49}{Let $\tau:I\to I$ be \pwm\ and continuous, with associated partition
$C$, and associated local homeo\-morphism $\sigma:X\to X$.   Let $B
=
\tau C
\setminus
\{0,1\}$.  Assume
\begin{enumerate}
\item The orbits of points in $B$ are infinite and
disjoint.\label{1.49i}
\item $\tau(\{0,1\}) \subset B$.\label{1.49ii}
\end{enumerate}
\noindent Then
 the only solution in $C(X,\Z)$ of
$\L g = g$ is $g = 0$.}

\prooff{We start by establishing some basic facts about $B$ that will also be used in the proof of Lemma \ref{1.50}. By
the definition of
$B$,
\begin{equation}
\tau C \subset B \cup \{0,1\},\label{(1.17)}
\end{equation}
and since $\tau(\{0,1\}) \subset B$,
\begin{equation}
\tau^2 C \subset \tau B \cup B.\label{(1.18)}
\end{equation}
Combining (\ref{(1.17)}) and (\ref{(1.18)}) gives
\begin{equation}
\orbit(C) \setminus C \subset \orbit(B).\label{(1.19)}
\end{equation}
By (\ref{(1.18)}) and the IDOC,
\begin{equation}
\tau^2(\orbit(B)) \cap \tau^2 C \,\subset\, \left(\bigcup_2^\infty \tau^n B\right) \cap (\tau B \cup B) = \emptyset.
\label{(1.19.1)}\end{equation}
Hence
\begin{equation}
\orbit(B) \cap C = \tau(\orbit B) \cap \tau C = \emptyset,\label{(1.20)} 
\end{equation}
since if (\ref{(1.20)}) failed, the result of applying $\tau^2$ or $\tau$ to each equality would contradict (\ref{(1.19.1)}).

By the
definition of
$I_1$, for each
$a
\in I_1$ there exists
$m
\ge 0$ such that $\tau^m a \in \orbit(C)$. By (\ref{(1.18)}), 
$\tau^2 C\subset
\orbit(B)$, and thus $\tau^{m+2}a\in \orbit(B)$.   Hence
\begin{equation}
a \in I_1 \implies  \tau^n a
\in\orbit(B)\hbox{ for some $n \ge 0$.}\label{(1.21)}
\end{equation}
\smallskip

Now assume $g \in C(X,\Z)$ and $\L g = g$. Choose $n\ge 0$ such that 
\begin{equation}
\tau^n \D(g) \subset \orbit(B).
\end{equation}
By (\ref{(1.13)}), (\ref{(1.17)}), and (\ref{(1.18)}), using $\L g  = g$,
\begin{equation}
\D(g) = \D(\L^n g) \,\subset \,\tau^n \D(g) \cup \bigcup_{k=1}^n \tau^k C
\,\subset\, \orbit(B) \cup \{0,1\}.\label{(1.22)}
\end{equation}

Suppose that  $\D(g) \not\subset \{0,1\}$, and let $N\ge 0$ be the largest integer
such that
$\D(g)$ contains a point of $\tau^N B$. (By the IDOC for $B$, the sets $\{\tau^n B
\mid n
\ge 0\}$ are disjoint, and $\D(g)$ is finite, so such an integer $N$ exists.) 
Choose $e \in \D(g) \cap \tau^NB$.
We now will verify that the hypotheses of Lemma \ref{1.45} hold here. Since $e \in \orbit(B)$, then $\tau e \notin \tau C$ by
(\ref{(1.20)}).
Let $x \in \D(g) \cap \tau^{-1}\tau e$. Then $\tau x = \tau e \notin \tau C$, so $x \notin C$. 
Thus $x \in \D(g)\setminus C$ and
(\ref{(1.22)}) imply that $x$ is in $\orbit (B)$.
 By the IDOC,
$\tau$ is 1-1 on $\orbit(B)$, so $\tau x= \tau e $  implies $x = e$. Now Lemma \ref{1.45}
implies that $\tau e \in \D(\L g)$. Then $\D(g)= \D(\L g)$ contains $\tau e \in \tau^{N+1}B$, which
contradicts the maximality of $N$.

We conclude that $\D(g) \subset \{0,1\}$, and therefore $g \in \Z I(0,1)$. If
$\L I(0,1) = I(0,1)$, since $\L I(0,1)(t)$  is the cardinality of $\sigma^{-1}(t)$,
then $\sigma$ would be bijective. It would follow that $\tau$ was bijective, and thus would be a homeo\-morphism, so
$\tau(\{0,1\})
\subset
\{0,1\}$, contrary to (\ref{1.49ii}). Thus $g = 0$.}

\lem{1.50}{Let $\tau:I\to I$ be \pwm\ and continuous, with associated partition
$C$, and associated local homeo\-morphism $\sigma:X\to X$.   Let 
$B=\tau C\setminus\{0,1\}$.  Assume
\begin{enumerate}
\item the orbits of points in $B$ are infinite and 
disjoint, and
\item $\tau(\{0,1\}) \subset B$.
\end{enumerate}
\noindent  If $f \in C(X,\Z)$ and $\D(f-\L f) \subset C$, then
$\D(f)\subset C$ and $\L f \in \Z I(0,1)$.}

\prooff{Assume $\D(f-\L f) \subset C$. We first prove
\begin{equation}
\D(f) \subset \orbit(C).\label{(1.23)}
\end{equation}
Suppose that $e_0 \in \D(f) \setminus \orbit(C)$. Observe that $e_0$ can't be
periodic, since by (\ref{(1.21)}), the orbit of $e_0$ lands in the orbit of $B$, and by
the IDOC, no point of $B$ can be periodic. Therefore the sets
$\{\tau^{-k}e_0\mid k \ge 0\}$ are disjoint, so there exists a maximal integer $N$ such that $\tau^{-N}e_0 \cap \D(f) \not=
\emptyset$. Choose $e  \in \tau^{-N}e_0 \cap \D(f)$, and observe that 
\begin{equation}
e \in \D(f) \setminus \orbit(C) \text{ and } \tau^{-1}e \cap \D(f) = \emptyset.\label{(1.24)}
\end{equation}
Then
$e
\notin \tau \D(f)$, and since $e \notin \orbit(C)$, in particular $e \notin \tau C$. By (\ref{(1.13)}),
\begin{equation}
\D(\L f) \subset \tau(\D(f)) \cup \tau C,\label{(1.24.1)}
\end{equation}
so $e \notin \D(\L f)$. Since $f$ is discontinuous at $e$ and $\L f$ is continuous at $e$, then $e$
is in $\D(f-\L f)$, which is contained in $C$ by hypothesis, contradicting $e \notin \orbit(C)$. This proves (\ref{(1.23)}).

We now prove 
\begin{equation}
\tau(\D(f) \setminus C) \subset  \D(f) \setminus C.\label{(1.25)}
\end{equation}
Let $e$ be an arbitrary member of $\D(f) \setminus C$. By (\ref{(1.23)}) and (\ref{(1.19)}),
\begin{equation}
e \in\D(f)\setminus C \subset \orbit(C) \setminus C \subset \orbit(B).
\end{equation}
By the same argument as in the proof of Lemma \ref{1.49}, the hypotheses of  Lemma \ref{1.45} are satisfied here, so $\tau e \in
\D(\L f)$. By (\ref{(1.20)}), $\tau e \in \orbit(B)$ implies that $\tau e \notin C$. By hypothesis, $\D(f-\L f) \subset C$, so
$\tau e \notin \D(f-\L f)$, i.e., $f-\L f$ is continuous at $\tau e$. Since $\tau e \in \D(\L f)$, then $\L f$ is discontinuous
at $\tau e$, so $f = (f-\L f) + \L f$ is discontinuous at $\tau e$. Hence $\tau e \in \D(f) \setminus C$,  which finishes the
proof of (\ref{(1.25)}).

 Since $\D(f)\setminus C \subset \orbit(B)$, by the
IDOC any point in
$\D(f) \setminus C$ has an infinite orbit, and this orbit stays in $\D(f)\setminus C$ by (\ref{(1.25)}). Since $\D(f)$ is
finite, it follows that   
$\D(f)\setminus C$ is empty, so we've proven
$$\D(f) \subset C.$$ 
By hypothesis $\D(f-\L f) \subset C$, so $\D(\L f) \subset C$. On the other hand, by (\ref{(1.24.1)}), 
\begin{equation}
\D(\L f) \subset \tau \D(f) \cup \tau C \subset \tau C.
\end{equation}
By (\ref{(1.20)}), $B\cap C= \emptyset$, so by (\ref{(1.17)}),
\begin{equation}
\D(\L f) \subset \tau C \cap C \subset (B \cup
\{0,1\})\cap C =
\{0,1\}.
\end{equation}
If follows that $\L f \subset \Z I(0,1)$, which completes the proof of the lemma.}

\lem{1.51}{Let $E_1, \ldots, E_n$ be elements in a $\Z[t]$ module $M$. If
\begin{enumerate}
\item $E_1, \ldots, E_n$ are independent in $M$ viewed as a  $\Z$ module,\label{1.51i}
\item $(1-t)M\cap \sum_i \Z E_i = \{0\}$, and\label{1.51ii}
\item $\ker(1-t) = \{0\}$,\label{1.51iii}
\end{enumerate}
\noindent then $E_1, E_2, \ldots, E_n$ are independent in $M$ viewed as a  $\Z[t]$ module.}

\prooff{We first prove the implication (for $p_1, p_2, \ldots, p_n \in \Z[t]$)
\begin{equation}
\sum_i p_i(t) E_i = 0 \implies p_1(1) = p_2(1) = \cdots = p_n(1) = 0.\label{(1.26)}
\end{equation}
If
\begin{equation}
0 = \sum_i p_i(t) E_i = \sum_i (p_i(t)-p_i(1))E_i + \sum_i p_i(1) E_i,
\end{equation}
since each $p_i(t) - p_i(1)$ is divisible by $(1-t)$, then
\begin{equation}
\sum_i p_i(1) E_i = - \sum_i (p_i(t)-p_i(1))E_i  \in (1-t)M,
\end{equation}
which by (\ref{1.51i}) and (\ref{1.51ii}) gives the implication (\ref{(1.26)}).

Now let $p_1, p_2, \ldots, p_n$ be arbitrary polynomials in $\Z[t]$, and assume that $\sum_i p_i(t)
E_i = 0$, and that not all $p_i$ are the zero polynomial. Let $k$ be the largest power such that
$(1-t)^k$ divides each $p_i(t)-p_i(1)$ for $1 \le i \le n$. Write $p_i(t) -p_i(1) = (1-t)^k
q_i(t)$, where $q_1, \ldots, q_n \in \Z[t]$. Note, for use below, that by the maximality of $k$,
there is an index
$j$ such that $(1-t)$ does not divide $q_j(t)$, and thus $q_j(1) \not= 0$.

Now by (\ref{1.51iii}),
\begin{equation}
0 = (1-t)^k \sum_i q_i(t)E_i
\end{equation}
implies that $\sum_i q_i(t) E_i = 0$, so by (\ref{(1.26)}), $q_1(1)= q_2(1) = \ldots = q_n(1) = 0$. This
contradicts $q_j(1) \not= 0$. Thus we conclude that $p_1 = p_2 = \cdots = p_n  = 0$, and hence that $E_1,
\ldots, E_n$ are independent in the $\Z[t]$ module $M$.}

\lem{1.52}{Let $\tau:I\to I$ be continuous, surjective, and \pwm, with associated partition
$C= \{a_0, a_1, \ldots, a_q\}$, and associated local homeo\-morphism $\sigma:X\to X$.  Assume that the orbits
of $a_1, a_2, \ldots, a_{q-1}$ are infinite and disjoint, and that $\tau(\{0,1\}) \cap \{0,1\} =
\emptyset$. Then $I(a_0,a_1)$, $I(a_1,a_2)$, $\ldots$, $I(a_{q-2},a_{q-1})$ are independent in
$DG(\tau)$.}

\prooff{Let $B = \tau C \setminus \{0,1\}$. By surjectivity,
\begin{equation}
\{0,1\} \subset \tau(\{a_1, \ldots, a_{q-1}\}),\label{(1.27)}
\end{equation}
and it follows that $B$ satisfies the
hypotheses of Lemma \ref{1.49} and \ref{1.50}.  Let $E_i = I(a_{i-1},a_i)$ for $1\le i \le q-1$. 

 We now
verify the hypotheses of Lemma \ref{1.51} with
$n = q-1$ and $M = C(X,\Z)$, viewed as a $\Z[t]$ module via the action of $\L = \L_\sigma$.

Since $E_1, \ldots, E_{q-1}$ are disjoint, 
if $\sum_i z_i E_i = 0$, where $z_1, \ldots, z_{q-1} \in \Z$, then $z_1 = \cdots = z_{q-1} =
0$, so $E_1,
\ldots, E_{q-1}$ are independent in $M$ viewed as a $\Z$ module.  

Next we show $(\id-\L)C(X,\Z) \cap \sum_i \Z E_i = \{0\}$. Suppose that $g \in  \sum_{i=1}^{q-1} \Z E_i$, with
$g\not= 0$, and that
$g = (\id-\L)f$ for some
$f \in C(X,\Z)$. Then $\D(f-\L f) = \D(g) \subset C$,  so by Lemma \ref{1.50}, $\D(f) \subset C$, and $\L f =z
I(0,1)$ for some
$z\in
\Z$.  Let
$f =
\sum_{i=1}^q z_i I(a_{i-1},a_i)$, where $z_1, \ldots, z_{q} \in \Z$. Then
\begin{equation}
z I(0,1) = \L f = \sum_{i=1}^q z_i I(\tau a_{i-1},\tau a_i).\label{(1.28)}
\end{equation}
By (\ref{(1.27)}) and the IDOC for $\{a_1, \ldots, a_{q-1}\}$, the points
$\{\tau a_j \mid 0 \le j \le q\}$ are all distinct. If $x$ is the first index in (\ref{(1.28)}) for which $z_x$ is
nonzero, and $y$ is the last such index, then the right side of (\ref{(1.28)}) will be discontinuous at
$\tau a_{x-1}$ and $\tau a_y$, so $\{\tau a_{x-1}, \tau a_y\} \subset \{0,1\}$, and $z_x = z_y = z\not= 0$.
Since
$g \in  \sum_{i=1}^{q-1} \Z E_i$ is continuous at $1$,
\begin{equation}
g = f -\L f = \sum_{i=x}^y z_i I(a_{i-1},a_i)-zI(0,1)
\end{equation}
implies that $a_y=1$. Since $\{\tau a_{x-1}, \tau a_y\} \subset \{0,1\}$, this
contradicts our assumption that $\tau(\{0,1\}) \cap \{0,1\} = \emptyset$, so we have shown that $g$
cannot be in $(\id-\L)C(X,\Z)$.

Finally, by Lemma \ref{1.49}, $\ker(\id-\L) =
\{0\}$, so we have verified all the assumptions of Lemma \ref{1.51}.  We conclude that $E_1,
\ldots, E_{q-1}$ are independent in the $\Z[t]$ module $C(X,\Z)$. It follows that their equivalence classes
are also independent in the
$\Z[t,t^{-1}] $ module $DG(\tau)$.}

\prop{1.53}{Let $\tau:I \to I$ be multimodal and surjective, with associated partition $C = \{a_0, a_1, a_2, \ldots,
a_q\}$.  If  the orbits of
$a_1,
\ldots, a_{q-1}$ are disjoint and infinite, and   $\tau(\{0,1\}) \cap \{0,1\}= \emptyset$,
then
$DG(\tau) \cong (\Z[t,t^{-1}])^{q-1}$ as (unordered) $\Z[t,t^{-1}]$ modules, with basis $I(a_0,a_1)$,
$I(a_1,a_2)$, $\ldots$, $I(a_{q-2},a_{q-1})$.}

\prooff{By Corollary \ref{1.42}, $I(a_0,a_1), I(a_1,a_2)\ldots, I(a_{q-1},a_q)$ generate $DG(\tau)$. We will show
$I(a_{q-1},a_q)$ is in the submodule generated by the first $q-1$ of these elements. It will follow that 
$I(a_0,a_1), I(a_1,a_2)\ldots, I(a_{q-2},a_{q-1})$ generate $DG(\tau)$, and are independent by Lemma
\ref{1.52}, and thus these elements form  a basis.

Since $\tau$ is surjective,
we can choose indices $x,y$ such that $\tau(a_x) = 0$ and
$\tau(a_y) = 1$.  For simplicity of
notation, we suppose that $x < y$.  Since  $\tau(\{0,1\}) \cap \{0,1\}= \emptyset$, then $0 < x < y < q$. 
Next observe that for $a,b \in I_1$, $I(0,a) \pm I(a,b) = I(0,b)$, with the $\pm$ sign depending on
whether $a < b$ or $a > b$.  Let $b_i = \tau(a_{x+i})$ for $0 \le i \le y-x$. In particular, $b_0 = \tau(a_x) = 0$ and
$b_{y-x} = \tau(a_y) = 1$. Then for appropriate choices of the $\pm$ signs in the equation below,
\begin{align}
{}&\L\left(I(a_x,a_{x+1}) \pm I(a_{x+1},a_{x+2})\pm\cdots\pm I(a_{y-1},a_y) \right)\notag \\
{}& =I(b_0,b_1) \pm I(b_1, b_2) \pm\cdots\pm I(b_{y-x-1},b_{y-x})\notag\\
{}& = I(b_0,b_{y-x}) = I(0,1).\notag
\end{align}
It follows that $I(0,1)$ is in the $\Z[t,t^{-1}]$ submodule of $DG(\tau)$ generated by $I(a_0,a_1)$,
$I(a_1,a_2)$, $\ldots$, $I(a_{q-2},a_{q-1})$. Since 
\begin{equation}
I(a_{q-1},a_q) = I(0,1) - I(a_0,a_1)-
I(a_1,a_2)-\cdots- I(a_{q-2},a_{q-1}),
\end{equation}
it follows that $I(a_0,a_1),
I(a_1,a_2)\ldots, I(a_{q-2},a_{q-1})$ generate $DG(\tau)$. This completes the proof that the prescribed
set of $q-1$ elements form a basis for the  module $DG(\tau)$, and thus that $DG(\tau) \cong (\Z[t,t^{-1}])^{q-1}$.}

\section{Interval exchange maps}
  
\defi{1.99}{A  \pwm\ map $\tau:[0,1)\to [0,1)$ is an {\it interval exchange map\/} if $\tau$
is linear with slope 1 on each interval of monotonicity, is  bijective, and is right
continuous at all points.  (See Figure \ref{fig1} on page~\pageref{fig1} for an example.) We
will usually identify $\tau$ with its extension to a map from
$[0,1]$ into
$[0,1]$, defined to be left continuous at 1.}

If $\{a_0, a_1, \ldots, a_n\}$ is the associated partition for an interval exchange map, then for each
$i$,  the lengths of $[a_{i-1},a_i)$ and $\tau([a_{i-1},a_i))$ will be the same, and the
images of the intervals
$\{\tau[a_{i-1},a_i)\mid 1\le i\le n\}$ will form a partition of $[0,1)$, which is the
reason for the name ``interval exchange map". 

\defi{1.62}{A map $\tau:I\to I$ is {\it essentially  injective\/} if there  are
no disjoint intervals $J, J'$ with $\tau(J) = \tau(J')$.}

This is equivalent to $\tau$ being injective on the complement of the set of  endpoints of intervals of
monotonicity.

\lem{1.63}{A \pwm\ map $\tau:I\to I$ is essentially injective iff the associated
local homeo\-morphism $\sigma:X\to X$ is injective.}

\prooff{If $J, J'$ are disjoint intervals in $I$ such that $\tau(J) = \tau(J')$, then
$\tau(J\cap I_0) = \tau(J'\cap I_0)$. Since $\sigma$ restricted to $X_0$  is conjugate to $\tau$ restricted to $I_0$, then
$\sigma$ is not 1-1. 

 Conversely, suppose
$\sigma$ is  not 1-1. Note that
for
$z
\in X$,
$(\L 1)(z)$ is the number of pre-images of $z$. Let $z\in X$ have at least two distinct pre-images. Since $\L
1$ is integer valued and continuous, then $(\L 1)(x) \ge 2$ for all $x$ in an open neighborhood $V$ of $z$. Thus
each member of $V$ has at least two pre-images.  It follows that each member of an open set
in
$X_0$ has two pre-images, so the same is true for $\tau$ on an open set in $I_0$. It
follows that
$\tau$ is not essentially injective.}

\defi{1.100}{A \pwm\ map $\tau:I\to I$ is {\it essentially bijective\/} if $\tau$ is surjective and
essentially injective, cf. Definitions \ref{1.23} and \ref{1.62}.}

 A \pwm\ map
$\tau:I\to I$ will be essentially bijective iff the associated local homeo\-morphism
$\sigma:X\to X$ is bijective, cf. Lemma \ref{1.63} and the remarks after Definition \ref{1.23}. Note that an interval
exchange map $\tau:I\to I$ is bijective on $[0,1)$, but might only be essentially bijective
on $I = [0,1]$.

 For our purposes, the fact that an interval exchange map is linear on
each interval of monotonicity will not play a role, which motivates the following definition.

\defi{1.101}{A  \pwm\ map $\tau:[0,1)\to [0,1)$ is a \emph{generalized interval exchange map} if
$\tau$ is increasing on each interval of monotonicity, is  bijective, and is right continuous at
all points.  We will usually identify $\tau$ with its extension to a map from $[0,1]$ into $[0,1]$, defined
to be left continuous at 1. (That map will be essentially bijective.)}

We are now going to find the dimension group for (generalized) interval exchange maps that have certain
orbits infinite and disjoint.  (See also Theorem \ref{1.88} and Proposition \ref{1.53} for continuous maps with the IDOC
property.)

\lem{1.102}{Let $\tau:I\to I$ be \pwm, with associated partition
$C= \{a_0, a_1, \ldots, a_n\}$, and associated local homeo\-morphism $\sigma:X\to X$. Let $E_i =
I(a_{i-1},a_i)$ for $1 \le i \le n-1$. If $\L= \L_\sigma$ is injective, and $\orbit(\htau C)
\cap  C
\subset
\{0,1\}$, then
$E_1,
\ldots, E_{n-1}$ are independent in $DG(\tau)$.}

\prooff{Suppose first, to reach a contradiction, that  
$E_1,
\ldots, E_{n-1}$ are not independent in the $\Z[t]$ module $C(X,\Z)$, and choose $p_1, \ldots, p_{n-1}$ in
$\Z[t]$, not all zero, satisfying
$\sum_i p_i(\L) E_i = 0$, with the sum  of the degrees of the non-zero
polynomials among $p_1,
\ldots, p_{n-1}$  as small as possible.  For $1 \le i \le n-1$, define $q_i(t)\in \Z[t]$ by
$t q_i(t) = p_i(t)-p_i(0)$, and define
\begin{equation}
g_1 = \sum_{i=1}^{n-1} p_i(0) E_i, \hbox{ and } g_2 = \sum_{i=1}^{n-1} q_i(\L) E_i.
\end{equation}
Then $g_1 + \L g_2 = \sum_i
p_i(\L) E_i = 0$.  Note by (\ref{(1.13)}),
\begin{equation}
\D(\L g_2) \subset \htau(\D(g_2)) \cup \htau C.
\end{equation}
Since $\cup_i \D(E_i) \subset C$, then $\D(g_2) \subset \orbit(C)$, so 
$\D(\L g_2) \subset
\orbit(\htau C)$. Since  $\D(g_1) \subset C\setminus \{1\}$ and $g_1 = -\L g_2$, by the assumption that
$\orbit(\htau C) \cap C \subset \{0,1\}$,  
\begin{equation}
\D(g_1) = \D(\L g_2) \subset \orbit(\htau C) \cap (C\setminus \{1\}) \subset \{0\}.
\end{equation}
Thus $\D(g_1) = \D(\L g_2) = \{0\}$, hence $g_1 = 0$, and
$\L g_2=\L
\sum_i q_i(\L)E_i = 0$. By assumption, $\L$ is injective, so $\sum_i q_i(\L)E_i = 0$, which contradicts the
minimality of the sum of the degrees of the non-zero polynomials among  $p_1, \ldots, p_{n-1}$. We conclude
that
$p_1,
\ldots, p_{n-1}$ are all the zero polynomial. Thus $E_1, \ldots, E_{n-1}$ are independent in the $\Z[t]$
module
$C(X,\Z)$. It follows that their equivalence classes are independent in the $\Z[t,t^{-1}]$ module
$DG(\tau)$.}

We will write $\tau(a^-)$ for $\lim_{t\to a^-} \tau(t)$, and define $\tau(a^+)$ similarly.

\prop{1.103}{Assume $\tau:I\to I$ is \pwm\ and essentially bijective, with
associated partition  $C= \{a_0, a_1, \ldots, a_n\}$.  If $\orbit(\htau C) \cap C \subset
\{0,1\}$, then the
(unordered) module
$DG(\tau)$ is isomorphic to $\Z[t,t^{-1}]^{n-1}\oplus
\Z$, where the action of
$\Z[t,t^{-1}]$ on the first summand is the usual one, and $t$ acts trivially on the
second summand.}

\prooff{Since $\tau$ is essentially bijective, the associated local homeo\-morphism
$\sigma:X\to X$ is bijective. Then for $g \in C(X,\Z)$, $\L_\sigma g = g\circ \sigma^{-1}$, so
$\L_\sigma$ is  bijective.

Let
 $E_i = I(a_{i-1},a_i)$ for $1 \le i \le n$. By Lemma \ref{1.102}, $E_1, \ldots, E_{n-1}$ are
independent in
$DG(\tau)$.
Let $\htau C = \{b_0, b_1, \ldots, b_n\}$, in increasing order. Then each interval $I(b_{j-1},b_j)$ is equal to $\L E_k$
for some $k$. The left and right limits of $\tau$ at partition points are in $\htau C$, so each
jump $I(\tau(a_i^-), \tau(a_i^+))$ at partition points is of the form $I(b_j,b_k)$ for some $j <
k$, and therefore is a sum of intervals of the form $\L E$ for $E \in \{E_1, \ldots, E_n\}$.
Furthermore, by surjectivity, $0 \in
\htau C$, so by Theorem \ref{1.41},  $\{E_1, \ldots, E_n\}$  is a set of generators for $DG(\tau)$.
Since $E_1+ E_2\ldots+ E_n=I(0,1)$, then the intervals $I(0,1)$, $E_1$, $E_2$, $\ldots$, $E_{n-1}$
also generate.

Since $\sigma$ is bijective, then $\L_*I(0,1) = I(0,1)$, so for every $p \in \Z[t,t^{-1}]$,
 $p(\L_*)I(0,1) = p(1)I(0,1)$. Thus $\Z[t,t^{-1}]I(0,1) = \Z I(0,1)$, so every element of $DG(\tau)$
 can be written in the
form $\sum_{i=1}^{n-1} p_i (\L_*) E_i + z I(0,1)$, with each $p_i \in \Z[t,t^{-1}]$, and $z\in \Z$. By the
independence of
$\{E_i \mid 1\le i \le n-1\}$, no nonzero element
$\sum_{i=1}^{n-1} p_i (\L_*) E_i$ is fixed by $\L_*$, so this representation is unique.}

When we have referred to the orbit of a point with respect to a \pwm\ map with some discontinuities, we
have previously used the orbit with respect to the multivalued map $\htau$. In the next proposition, we
need to refer to the orbit with respect to $\tau$ itself, and we will call this the $\tau$-orbit to avoid
confusion with our earlier usage.  Recall that if a set $B$ has the property that its points have orbits
that are infinite and disjoint, we refer to this property as the IDOC for $B$.

\cor{1.104}{Let $\tau:I\to I$ be a generalized interval exchange map, with associated partition
$C= \{a_0, a_1, \ldots, a_n\}$.  If the $\tau$-orbits of $a_1, a_2, \ldots,
a_{n-1}$ are infinite and disjoint, then $DG(\tau)= C(X,\Z)$ is isomorphic as an (unordered) $\Z[t,t^{-1}]$
module to
$\Z[t,t^{-1}]^{n-1}\oplus
\Z$ (where the action of $\L_*$ on the first summand is coordinate-wise multiplication by $t$, and on the second
summand is trivial).}

\prooff{Let $C_0 = C\setminus \{0,1\} = \{a_1, a_2, \ldots, a_{n-1}\}$. We'll prove
\begin{equation}\label{eq104} x \in C \implies \orbit(\htau x)
\cap C_0 = \emptyset.
\end{equation}
First consider the case $x = 0$. If $\tau 0 = 0$, then (\ref{eq104}) is clear. If $\tau 0 \not= 0$, by surjectivity, we can
choose
$a \in C_0$ with $\tau a = 0$. Then $\orbit(\htau 0 ) = \orbit(\tau^2 a)$, and $\orbit(\tau^2 a)$ misses $C_0$ by the IDOC.

Now suppose $x \in C_0$. Then the orbit of $\tau x$ misses $C_0$ by the IDOC. Furthermore, if $x \in C_0 \cup \{1\}$, there
exists $y \in C_0 \cup \{0\}$ such that $\tau y = \tau(x^-)$. Then $\orbit(\tau(x^-)) = \orbit(\tau y)$ misses $C_0$ by the
previous part of this proof. Thus we've shown $\orbit(\htau C) \cap C_0 = \emptyset$.

As observed in the proof of Proposition \ref{1.103}, if $\sigma:X\to X$ is the local homeo\-morphism
associated with $\tau$, then essential bijectivity of
$\tau$ implies that
$\L=
\L_\sigma$ is bijective, and so
$DG(\tau) = C(X,\Z)$. The rest of the corollary follows from Proposition \ref{1.103}.}

Note that in Corollary \ref{1.104}, the order on $DG(\tau)= C(X,\Z)$ is just the usual pointwise order on $C(X,\Z)$.

\section{C*-algebras}

Let $\tau:I\to I$ be a \pwm\ map, with associated partition $C= \{a_0, a_1, \ldots, a_n\}$. We now
give a dynamical construction of  a C*-algebra
$A_\tau$ such that
$K_0(A_\tau)
\cong DG(\tau)$. The construction will also display $DG(\tau)$ explicitly as an inductive limit of groups of the form
$\Z^{n_k}$, with the connecting maps easily computable from the map $\tau$.

 Define $C_0 = \{0,1\}$, and for
$n \ge 1$, define
$C_n =
(\bigcup_{k=1}^{2n} \htau^k C) \bigcap \htau^n(I) $.  Observe that the sequence $\{C_n\}$ has the following
properties.
\begin{align}
\htau(C) &\subset C_1.\label{(1.38)}\\
\htau(C_n) &\subset C_{n+1}.\label{(1.39)}\\
C_n&\subset \htau^n I.\label{(1.39.1)}\\
C_k \cap \htau^{n}I &\subset C_{n} \text{ for $k\le n$}.\label{(1.39.2)}
\end{align}
\vskip-.78cm
\begin{equation}
\hbox{If $a \in I_1$, then there exists $n \ge 0$ such that $\htau^n(a) \subset C_n$.}\label{(1.40)}
\end{equation}
Any sequence
$\{C_n\}$ of finite subsets of
$I_1$ satisfying these five properties would work equally well in what follows. (Note that
(\ref{(1.40)}) could fail if we replaced the upper limit $2n$ by $n$ in the union defining $C_n$.) 
We say $x,y \in
C_n$ are {\it adjacent in $C_n$\/} if there is no point in $C_n$ strictly between
$x$ and $y$. For $n \ge 0$, let
$\P_n$ be the set of closed intervals contained in $\htau^n(I)$
whose endpoints are adjacent points in $C_n$. For each
$Y
\in
\P_n$, choose a point
$x$ in the interior of $Y$ and define $k(Y)$ to be the
cardinality of $\tau^{-n}(x)$.  (By the definition of $\P_n$, $k(Y)$
does not depend on the choice of $x$.)  For $n \ge 0$, define
\begin{equation}
A_n = \oplus_{Y\in \P_n} C(Y,M_{k(Y)}).
\end{equation}

Recall that for each $i$,  $\tau_i$ denotes the continuous extension  of $\tau|_{(a_{i-1}, a_i)}$  to a
homeo\-morphism on $[a_{i-1},a_i]$.  Fix $n \ge 0$. If
$Z\in
\P_{n+1}$, for each index
$i$, 
$\tau_i(I)$ either contains $Z$ or is disjoint from the interior of $Z$. In the former case,
there is exactly one $Y \in \P_n$ such that such that
$\tau_i(Y) \supset Z$, or equivalently, $\tau_i^{-1}(Z) \subset Y$. 

For $n \ge 0$, let $X^n$ be the disjoint union of $\{Y \in \P_n\}$. We view $A_n$ as a subalgebra of
$C(X^n, \oplus_{Y\in \P_n} M_{k(Y)})$. Define $\phi_n: A_n\to A_{n+1}$ by
\begin{equation}
\phi_n(f) =\oplus_{Z\in \P_{n+1}}\hbox{diag }
(f\circ 
\tau_{i_1}^{-1}, \ldots, f\circ 
\tau_{i_p}^{-1}),\label{(1.41)}
\end{equation}
where for each $Z\in \P_{n+1}$ the indices $i_1, \ldots, i_p$ are the indices $i$ such that
$\tau_i(I)$ contains $Z$.

Let $A_\tau$ be the inductive limit of the sequence
$\phi_n:A_n\to A_{n+1}$.  Recall that an {\it interval algebra\/} is one isomorphic to $C([0,1], F)$,
where
$F$ is a finite dimensional C*-algebra, and an AI-algebra is an inductive limit of a sequence of
interval algebras. Thus each $A_n$ is an interval algebra, and $A_\tau$ is an AI-algebra.

\exe{1.111}{We illustrate this construction with the full tent map, cf. Example \ref{1.92}.  Here $C_n
=
\{0,1\}$ for all $n$, 
$\tau_1^{-1}(x) = x/2$ and
$\tau_2^{-1}(x) = 1-x/2$. Hence for all $n$, $A_n = C([0,1], M_{2^n})$, and $\phi_n:A_n\to A_{n+1}$ is given by
$\phi_n(f) = {\rm diag}(f\circ \tau_1^{-1}, f\circ \tau_2^{-1})$. (This example also appears in \cite{BratEll} and
\cite{Deac}.)}

\exe{1.112}{Let $\tau$ be the restricted tent map $\tau = T_{\sqrt{2}}$. 
In this case the critical point $c$ is eventually fixed, with $\tau^4 c = \tau^3 c$. Let $c= c_0, c_1, c_2, c_3,
\cdots 
$ be the orbit of $c$, and let $p$ be the fixed point of $\tau$. Note that $0=c_2 < c_0 = c < c_3 = p < c_1 = 1$.  Then we can
represent the imbeddings by the following diagram of intervals.

\begin{figure}[htb]
\centerline{\includegraphics{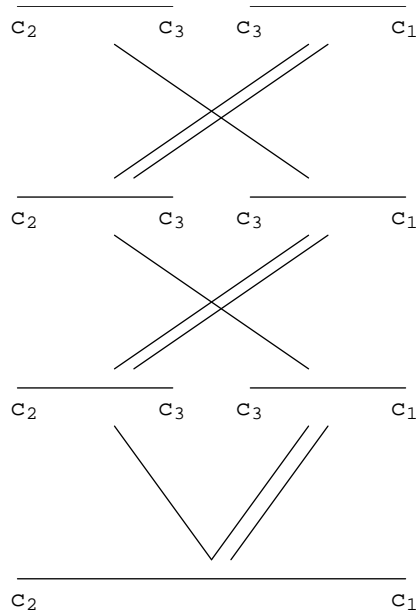}}
\caption{\label{fig2} Inductive limit construction}
\end{figure}

The disjoint union of the intervals on each level of this diagram are the sets
$X^0, X^1, X^2,
\ldots$ described previously, with the connecting lines indicating the maps of intervals by $\tau$. 
From this diagram, the algebras $A_0, A_1, \ldots$ and the maps $\phi_n:A_n\to A_{n+1}$ can also be
read off. For example, 
\begin{equation}
A_3 =  C([c_2,c_3],M_2) \oplus C([c_3,c_1],M_4).
\end{equation}
Here the dimensions of the matrix algebras for a given interval are the number of paths from the
bottom to that interval. 
Since $K_0(C([0,1], M_n)) = \Z$, $K_0(A_n)$ is found by  replacing each interval by a
copy of
$\Z$.  The
connecting lines indicate the connecting maps for the inductive limit of dimension groups, giving $DG(\tau)$.
 In this
particular example, the map $\tau$ is Markov, so the width of the diagram stabilizes. For general
$\tau$ the diagram can become steadily wider, as can be seen from the examples discussed in the previous
section.}

The diagram in Figure \ref{fig2}, with intervals shrunk to points,  also  can be thought of as a Bratteli
diagram, cf. \cite{Brat}. The  AF-algebra associated with that diagram would have the same $K_0$ group as
$A_\tau$.  However, $A_\tau$  is not necessarily an AF-algebra. 

\prop{1.110}{If $\tau:I\to I$ is \pwm, $K_0(A_\tau)$ and $DG(\tau)$ are isomorphic dimension groups.}

\prooff{Recall that for each $n$, we view elements of $A_n$ as maps from $X^n$ into a suitable matrix
algebra.  We view elements of the matrix algebra $M_q(A_n)$ as functions from $X^n$
into $\oplus_{Y\in \P_n} M_q(M_{k(Y)})$. For each projection
$p$ in a matrix algebra over $A_n$, define $\dim(p)$ to be the function $y\mapsto \dim\ p(y)$. Then
$\dim(p)
\in C(X^n,\Z)$, and since $X_n$ is a disjoint union of closed intervals,  $\dim$ extends uniquely to an isomorphism of
$K_0(A_n)$ onto
$C(X^n,\Z)$.  

Let
$\pi_n:\sigma^n(X)\to X^n$ be the collapse map (that collapses $x^\pm$ to the single point $x$ except for
the pre-images of points in
$C_n$.) For each $n\ge 0$, let $\widehat C(X^n,\Z)$ denote the subgroup of $C(X,\Z)$ generated by
characteristic functions of intervals $I(a,b)$ with $a, b \in C_n$. Then
$f\mapsto f\circ
\pi_n$ is an order isomorphism from $C(X^n,\Z)$  onto $\widehat C(X^n,\Z)$.  Observe that 
$$\widehat C(X^n,\Z) = \{f \in C(X,\Z) \mid \D(f) \subset C_n\},$$ 
where $\D(f)$ is the set of discontinuities of $f$ viewed as a function on $\R$, cf. Section \ref{Cyclic}. If $f \in  \widehat
C_n(X,\Z)$, then
$\supp f
\subset
\sigma^n(X)$, and
$\supp \L f \subset \sigma^{n+1}(X)$. By
(\ref{(1.13)}) and (\ref{(1.38)})--(\ref{(1.39.2)}),
$$\D(\L f) \subset \bigl(\htau(\D(f)) \cup \htau C\bigr)\cap \htau^{n+1}I  \subset (\htau(C_n) \cup \htau C)\cap
\htau^{n+1}I\subset C_{n+1}.$$ Thus
$\L$ maps $\widehat C_n(X,\Z)$ into $\widehat C_{n+1}(X,\Z)$.

 Consider
the following  commutative diagram.
\begin{equation}
\begin{matrix}K_0(A_n)& \mapright{\dim}& C(X^n,\Z) & \mapright{} &\widehat C_n(X,\Z)\cr
& &  &&\cr
\mapdown{(\phi_n)_*}&&&&\mapdown{\L}\cr
& &  &&\cr
K_0(A_{n+1}) & \mapright{\dim}& C(X^{n+1},\Z) & \mapright{} &\widehat C_{n+1}(X,\Z)
\end{matrix}
\end{equation}
(The commutativity of this diagram can  be verified by checking it on generators of $K_0(A_n)$, which are of the form $[1_Y
e_{11}]$ for $Y \in \P_n$.) Since the horizontal maps are isomorphisms, by continuity of
$K_0$ with respective to inductive limits,
$K_0(A)$ is isomorphic to the inductive limit of the sequence 
$\L:\widehat C_n(X,\Z)\to  \widehat C_{n+1}(X,\Z)$. We will be done if we show this inductive limit is
isomorphic to $DG(\tau)$. 

For that purpose, we apply Lemma \ref{1.81.1}. 
 By (\ref{(1.13)}), for $f \in C(X,\Z)$,
\begin{align} \label{(1.45)}
\begin{split}\D(\L^n f) &\subset \htau^n\D(f) \cup \biggl(\bigl(\bigcup_{k=1}^n \htau^k C \bigr)\cap \htau^n
I\biggr)\\
 & \subset \htau^n \D(f) \cup C_n. 
\end{split}
\end{align}
By (\ref{(1.39)}) and (\ref{(1.40)}), for each $f \in C(X,\Z)$ there exists $n \ge 0$
such that $\htau^n \D(f) \subset C_n$, so (\ref{(1.45)}) implies that $\L^n f \in
\widehat C_n(X,\Z)$. If $f \in \widehat C_k(X,\Z)$, then $\supp f \subset \sigma^k(X)$, so $f \in \L^k C(X,\Z)$ (Lemma
\ref{1.14}).   By Lemma
\ref{1.81.1},
$DG(\tau)$ is isomorphic to the inductive limit of the sequence
$\L:C_n(X,\Z)
\to C_{n+1}(X,\Z)$, so
$K_0(A)$ is isomorphic to
$DG(\tau)$.}

\end{document}